# GENERALIZED RIEMANN FUNCTIONS, THEIR WEIGHTS, AND THE COMPLETE GRAPH

NICOLAS FOLINSBEE AND JOEL FRIEDMAN

ABSTRACT. By a *Riemann function* we mean a function $f\colon \mathbb{Z}^n \to \mathbb{Z}$ such that $f(\mathbf{d})$ is equals 0 for $d_1+\cdots+d_n$ sufficiently small, and equals $d_1+\cdots+d_n+C$ for a constant, $C$, for $d_1+\cdots+d_n$ sufficiently large. By adding 1 to the Baker-Norine rank function of a graph, one gets an equivalent Riemann function, and similarly for related rank functions.

To each Riemann function we associate a related function $W\colon \mathbb{Z}^n \to \mathbb{Z}$ via Möbius inversion that we call the *weight* of the Riemann function. We give evidence that the weight seems to organize the structure of a Riemann function in a simpler way: first, a Riemann function $f$ satisfies a Riemann-Roch formula iff its weight satisfies a simpler symmetry condition. Second, we will calculate the weight of the Baker-Norine rank for certain graphs and show that the weight function is quite simple to describe; we do this for graphs on two vertices and for the complete graph.

For the complete graph, we build on the work of Cori and Le Borgne who gave a linear time method to compute the Baker-Norine rank of the complete graph. The associated weight function has a simple formula and is extremely sparse (i.e., mostly zero). Our computation of the weight function leads to another linear time algorithm to compute the Baker-Norine rank, via a formula likely related to one of Cori and Le Borgne, but seemingly simpler, namely

$$r_{\mathrm{BN},K_n}(\mathbf{d}) = -1 + \left| \left\{ i = 0,\ldots, \deg(\mathbf{d}) \ \Big| \ \sum_{j=1}^{n-2} \big( (d_j - d_{n-1} + i) \bmod n \big) \leq \deg(\mathbf{d}) - i \right\} \right|.$$

Our study of weight functions leads to a natural generalization of Riemann functions, with many of the same properties exhibited by Riemann functions.

## Contents

*Date*: Monday 30$^{\text{th}}$ May, 2022.
2010 *Mathematics Subject Classification.* Primary: 05C99.
Research supported in part by an NSERC grant.
Research supported in part by an NSERC grant.

1. INTRODUCTION

The main goal of this article is to give a combinatorial study of what we call *Riemann functions* and their *weights*. Our main motivation is to gain insight into the special case that is the Graph Riemann-Roch fomula of Baker and Norine [BN07]; the Baker-Norine formula has received a lot of recent attention [CB13, Bac17, MS14, CLM15], as has its generalization to *tropical curves* and other settings in recent years [Bac17, GK08, HKN13, JM13, AC13, MS13, AM10, CDPR12].

We were first interested in weights to address a question posed in [BN07] regarding whether or not their Graph Riemann-Roch formula could be understood as an Euler characteristic equation; this is partially answered in [FF]. However, weights are interesting for a number of purely combinatorial reasons: first, a Riemann-Roch formula is simpler to express in terms of the weight of the Riemann function. Second, the weights of the Riemann-Roch functions of certain graphs are very simple to write down. For example, in this article we build on the methods of Cori and Le Borgne [CB13] to give a very simple formula for the weights of the Baker-Norine rank function of a complete graph; this will allow us to prove a likely simpler variant of their algorithm to compute the values of this rank function. Furthermore, for the above reasons, as well as its connections to sheaves and Euler characteristics in [FF], we suspect that weights may be a useful way to describe many Riemann functions.

This article has two types of results: foundational results on Riemann functions and Riemann-Roch type formulas, and calculations of the weights of Baker-Norine rank functions of two types of graphs. Let us briefly summarize the results, assuming some terminology that will be made precise in Section 2.

1.1. **Riemann Functions and Weights.** By a *Riemann function* we mean a function $f\colon \mathbb{Z}^n \to \mathbb{Z}$ such that $f(\mathbf{d}) = f(d_1, \ldots, d_n)$ is *initially zero*, meaning $f(\mathbf{d}) = 0$ for $\deg(\mathbf{d}) = d_1 + \cdots + d_n$ sufficiently small, and *eventually*—meaning for $\deg(\mathbf{d})$ sufficiently large—equals $\deg(\mathbf{d}) + C$ for a constant, $C \in \mathbb{Z}$, which we call the *offset of $f$*. By adding 1 to the Baker-Norine rank function of a graph, one gets an equivalent Riemann function, and similarly for related rank functions.

If $f\colon \mathbb{Z}^n \to \mathbb{Z}$ is any function that is initially zero, then there is a unique, initially zero $W$ such that
$$f(\mathbf{d}) = \sum_{\mathbf{d}' \leq \mathbf{d}} W(\mathbf{d}')$$
where $\leq$ the usual partial order on $\mathbb{Z}^n$ (i.e., $\mathbf{d}' \leq \mathbf{d}$ means $d'_i \leq d_i$ for all $i = 1, \ldots, n$); we call $W$ the *weight* of $f$. If $f$ is a Riemann function, then $W$ is also eventually zero; much of what we prove about Riemann functions also holds for *generalized Riemann functions*, which we define as any initially zero function $f$ whose weight is eventually zero.

Returning to a Riemann function $f\colon \mathbb{Z}^n \to \mathbb{Z}$ with offset $C$, for any $\mathbf{K} \in \mathbb{Z}^n$ there exists a unique function $f_\mathbf{K}^\wedge$ such that for all $\mathbf{d} \in \mathbb{Z}^n$ we have

(1) $$f(\mathbf{d}) - f_\mathbf{K}^\wedge(\mathbf{K} - \mathbf{d}) = \deg(\mathbf{d}) + C,$$

and we refer to as a *generalized Riemann-Roch formula*; $f_\mathbf{K}^\wedge$ is also a Riemann function. Furthermore, if $f_\mathbf{K}^\wedge = f$ for some $f, K$, then the formula reads

$$f(\mathbf{d}) - f(\mathbf{K} - \mathbf{d}) = \deg(\mathbf{d}) + C,$$

which is the usual type of Riemann-Roch formula, both the classical formula of Riemann-Roch, and the Baker-Norine analog. Hence, our view of Riemann-Roch formulas is more "happy-go-lucky" than is common in the literature: for each $f, \mathbf{K}$ there is a generalized Riemann-Roch formula (1); we study any such formula, and view the case where $f_\mathbf{K}^\wedge = f$ as a special case which we call *self-duality*.

We are interested in weight functions, $W$, for a number of reasons:

(1) the weights of the Baker-Norine rank (plus 1) of the graphs we study in this article turn out be be simple to describe and very sparse (i.e., mostly 0); by contrast, at least for the complete graph, the Baker-Norine function is more difficult to compute. Hence the weights may be a more efficient way to encode certain Riemann functions of interest.
(2) For a Riemann function $f\colon \mathbb{Z}^n \to \mathbb{Z}$, the weight of $f_\mathbf{K}^\wedge$ turns out to equal $(-1)^n W_\mathbf{L}^*$, where $\mathbf{L} = \mathbf{K} + \mathbf{1}$ (where $\mathbf{1} = (1, \ldots, 1)$), and $W_\mathbf{L}^*$ is the function $W_{\mathbf{L}}^*(\mathbf{d}) = W(\mathbf{L} - \mathbf{d})$; hence it seems easier to check self-duality using the weight, $W$, rather than directly on $f$.

(3) In [FF], we model Riemann functions by restricting $f \colon \mathbb{Z}^n \to \mathbb{Z}$ to two of its variables, while holding the other $n-2$ variables fixed; if $f$ satisfies self-duality, a two-variable restriction, $\widetilde{f} \colon \mathbb{Z}^2 \to \mathbb{Z}$, of $f$ will generally not be self-dual; however $\widetilde{\mathbf{K}} \in \mathbb{Z}^2$ can be described as a restriction of $f_{\mathbf{K}}^{\wedge}$ (for any $\mathbf{K} \in \mathbb{Z}^n$). Since self-duality isn't preserved under restrictions, but generalized Riemann-Roch formulas behave well under restrictions, it seems essential to work with generalized Riemann-Roch formulas (1) in [FF] or whenever we wish to work with restrictions of Riemann functions to a subset of their variables.

(4) In certain Riemann functions of interest, such as those considered by Amini and Manjunath [AM10], self-duality does not generally hold, and yet one can always work with weights and generalized Riemann-Roch formulas.

(5) The formalism of weights applies to generalized Riemann functions, which is a much wider class of functions, and we believe likely to be useful in future work to model other interesting functions. In this case (1) is replaced by

$$f(\mathbf{d}) - f_{\mathbf{K}}^{\wedge}(\mathbf{K} - \mathbf{d}) = h(\mathbf{d}),$$

where $h$ is the unique *modular function* that eventually equals $f$ (see Section 3). One might expect such formulas to hold when, for example $f = f(\mathbf{d})$ is the sum of even Betti numbers of a sheaf depending on a parameter $\mathbf{d} \in \mathbb{Z}^n$, whose Euler characteristic equals a modular function $h$.

1.2. **The Weight of the Baker-Norine rank for Two Types of Graphs.** The second type of result in this article concerns the weights of the Baker-Norine rank function (plus 1) for two types of graphs, namely graphs on two vertices and the complete graph, $K_n$, on $n$ vertices. Both types of weight functions are quite simple and very sparse (i.e., mostly 0). For $K_n$ we build on the ideas of Cori and Le Borgne [CB13] to compute the weight of the Baker-Norine rank. A side effect of this computation is a formula for the Baker-Norine rank:

$$r_{\mathrm{BN},K_n}(\mathbf{d}) = -1 + \left| \left\{ i = 0, \ldots, \deg(\mathbf{d}) \,\bigg|\, \sum_{j=1}^{n-2} ((d_j - d_{n-1} + i) \bmod n) \leq \deg(\mathbf{d}) - i \right\} \right|,$$

where the "mod" function above returns a value in $\{0, \ldots, n-1\}$; this looks related to a formula given by Cori and Le Borgne. We also explain that—like the Cori and Le Borgne algorithm—there is an algorithm that computes this function in time $O(n)$. Our proof of this formula is self-contained, although uses some of the observations of Cori and Le Borge including one short and rather ingenious idea of theirs regarding the Baker-Norine function on a complete graph.

1.3. **Organization of this Article.** The rest of this article is organized as follows. In Section 2 we give some basic terminology, including the definition of a *Riemann function* and some examples, which (after subtracting 1) includes the Baker-Norine rank. In Section 3 we discuss what we mean by the *weight* of a Riemann function; this leads to a notation of *generalized Riemann functions*, which share many of the properties of Riemann functions. In Section 4 we define what we mean by a Riemann-Roch formula; we describe the equivalent condition on weights, which is simpler; these ideas generalize in a natural way to the setting of generalized

Riemann functions. In Section 5 we compute the weight of the Baker-Norine rank for graphs on two vertices, joined by any number of edges. In Section 6 we compute the weight of the Baker-Norine rank for a complete graph on $n$ vertices, and we give a formula for the Baker-Norine rank, which—like a related formula of Cori and Le Borgne—allows the rank to be computed in linear time in $n$. In Section 7 we prove our main theorems—stated earlier—that characterize *modular functions* used to define generalized Riemann functions.

## 2. Basic Terminology and Riemann Functions

In this section we introduce some basic terminology and define the notion of a Riemann function. Then we give some examples of Riemann functions.

### 2.1. Basic Notation.
We use $\mathbb{Z}, \mathbb{N}$ to denote the integers and positive integers; for $a \in \mathbb{Z}$, we use $\mathbb{Z}_{\le a}$ to denote the integers less than or equal to $a$, and similarly for the subscript $\ge a$. For $n \in \mathbb{N}$ we use $[n]$ to denote $\{1, \ldots, n\}$. We use bold face $\mathbf{d} = (d_1, \ldots, d_n)$ to denote elements of $\mathbb{Z}^n$, using plain face for the components of $\mathbf{d}$; by the *degree* of $\mathbf{d}$, denoted $\deg(\mathbf{d})$ or at times $|\mathbf{d}|$, we mean $d_1 + \ldots + d_n$.

We set
$$\mathbb{Z}^n_{\deg 0} = \{\mathbf{d} \in \mathbb{Z}^n \mid \deg(\mathbf{d}) = 0\},$$
and for $a \in \mathbb{Z}$ we similarly set
$$\mathbb{Z}^n_{\deg a} = \{\mathbf{d} \in \mathbb{Z}^n \mid \deg(\mathbf{d}) = a\}, \quad \mathbb{Z}^n_{\deg \le a} = \{\mathbf{d} \in \mathbb{Z}^n \mid \deg(\mathbf{d}) \le a\}.$$

We use $\mathbf{e}_i \in \mathbb{Z}^n$ (with $n$ understood) be the $i$-th standard basis vector (i.e., whose $j$-th component is 1 if $j = i$ and 0 otherwise), and for $I \subset [n]$ (with $n$ understood) we set
$$\mathbf{e}_I = \sum_{i \in I} \mathbf{e}_i; \tag{2}$$
hence in case $I = \emptyset$ is the empty set, then $\mathbf{e}_\emptyset = \mathbf{0} = (0, \ldots, 0)$, and similarly $e_{[n]} = \mathbf{1} = (1, \ldots, 1)$.

For $n \in \mathbb{N}$, we endow $\mathbb{Z}^n$ with the usual partial order, that is
$$\mathbf{d}' \le \mathbf{d} \quad \text{iff} \quad d'_i \le d_i \ \forall i \in [n],$$
where $[n] = \{1, 2, \ldots, n\}$.

### 2.2. Riemann Functions.
In this section we define *Riemann functions* and give examples that have appeared in the literature.

**Definition 2.1.** We say that a function $f \colon \mathbb{Z}^n \to \mathbb{Z}$ is a Riemann function if for some $C, a, b \in \mathbb{Z}$ we have
 (1) $f(\mathbf{d}) = 0$ if $\deg(\mathbf{d}) \le a$; and
 (2) $f(\mathbf{d}) = \deg(\mathbf{d}) + C$ if $\deg(\mathbf{d}) \ge b$;
we refer to $C$ as the *offset* of $f$.

In our study of Riemann functions, it will be useful to introduce the following terminology.

**Definition 2.2.** If $f, g$ are functions $\mathbb{Z}^n \to \mathbb{Z}$, we say that $f$ *equals $g$ initially* (respectively, *eventually*) if $f(\mathbf{d}) = g(\mathbf{d})$ for $\deg(\mathbf{d})$ sufficiently small (respectively, sufficiently large); similarly, we say that that $f$ is *initially zero* (respectively *eventually zero*) if $f(\mathbf{d}) = 0$ for $\deg(\mathbf{d})$ sufficiently small (respectively, sufficiently large).

Therefore $f\colon \mathbb{Z}^n \to \mathbb{Z}$ is a Riemann function iff it is initially zero and it eventually equals the function $\deg(\mathbf{d}) + C$, where $C$ is the offset of $f$.

2.3. **The Baker-Norine Rank and Riemann-Roch Formula.** In this article we study examples of the Baker-Norine rank for various graphs. In this subsection we briefly review its definition and its properties; for more details, see [BN07].

We will consider graphs, $G = (V, E)$ that are connected and may have multiple edges but no self-loops. Recall that if $G = (V, E)$ is any graph, then its *Laplacian*, $\Delta_G$ equals $D_G - A_G$ where $D_G$ is the diagonal degree counting matrix of $G$, and $A_G$ is the adjacency matrix of $G$.

**Definition 2.3** (The Baker-Norine rank function of a graph). Let $G = (V, E)$ be a connected graph without self-loops (but possibly multiple edges) on $n$ vertices that are ordered as $v_1, \ldots, v_n$. Hence we view its Laplacian, $\Delta_G$, as a map $\mathbb{Z}^n \to \mathbb{Z}^n$. Let $L = \mathrm{Image}(\Delta)$. We say that $\mathbf{d}, \mathbf{d}' \in \mathbb{Z}^n$ are *equivalent*, written $\mathbf{d} \sim \mathbf{d}'$, if $\mathbf{d} - \mathbf{d}' \in L$, and say that $\mathbf{d}$ is *effective* if $\mathbf{d} \geq \mathbf{0}$. Let $\mathcal{N}$ be the elements of $\mathbb{Z}^n$ that are not equivalent to an effective element of $\mathbb{Z}^n$; in particular
$$\deg(\mathbf{d}) < 0 \Rightarrow \mathbf{d} \in \mathcal{N}.$$

Consider
$$(3) \qquad f(\mathbf{d}) = \rho_{L^1}(\mathbf{d}, \mathcal{N}) = \min_{\mathbf{d}' \in \mathcal{N}} \|\mathbf{d} - \mathbf{d}'\|_{L^1},$$

where $\|\cdot\|_{L^1}$ is the usual $L^1$-norm
$$\|(x_1, \ldots, x_n)\|_{L^1} = |x_1| + \cdots + |x_n|.$$
We also write $f = f_G$, to emphasize the graph $G$, although its definition as a function $\mathbb{Z}^n \to \mathbb{Z}$ also depends on the ordering $v_1, \ldots, v_n$ of its vertices. The *Baker-Norine rank* of $\mathbf{d}$, denoted $r_{\mathrm{BN}}(\mathbf{d})$, is $f(\mathbf{d}) - 1$.

Since $f(\mathbf{d}) = 0$ iff $\mathbf{d} \in \mathcal{N}$, which is the case if $\deg(\mathbf{d}) < 0$, it follows $f$ is initially zero, and hence $r_{\mathrm{BN}}(\mathbf{d})$ initially equals $-1$. We remark that for $f(\mathbf{d}) \geq 0$ we easily see that both:

(1) $f(\mathbf{d})$ equals the largest integer $m \geq 0$ such that for any $\mathbf{a} \geq \mathbf{0}$ and of degree $m$ we have that $\mathbf{d} - \mathbf{a}$ is equivalent to an effective element of $\mathbb{Z}^n$, and
(2) $f(\mathbf{d}) = 1 + \min_{i \in [n]} f(\mathbf{d} - \mathbf{e}_i)$.

The Baker-Norine *Graph Riemann-Roch* formula states that for all $\mathbf{d}$ we have
$$(4) \qquad r_{\mathrm{BN}}(\mathbf{d}) - r_{\mathrm{BN}}(\mathbf{K} - \mathbf{d}) = \deg(\mathbf{d}) + 1 - g$$
where
(1) $g = 1 + |E| - |V|$ (which is non-negative since $G$ is connected), and
(2) $\mathbf{K} = \bigl(\deg_G(v_1) - 2, \ldots, \deg_G(v_n) - 2\bigr)$, where $\deg_G(v)$ is the degree of $v$ in $G$, i.e., the number of edges incident upon $v$ in $G$.

It follows that for all $\mathbf{d} \in \mathbb{Z}^n$
$$(5) \qquad f(\mathbf{d}) - f(\mathbf{K} - \mathbf{d}) = \deg(\mathbf{d}) + 1 - g.$$
It follows that for $\mathbf{d}$ such that
$$\deg(\mathbf{d}) > \deg(\mathbf{K}) = \sum_i \Bigl(\deg_G(v_i) - 2\Bigr) = 2|E| - 2|V|$$
we have $f(\mathbf{K} - \mathbf{d}) = 0$; hence
$$(6) \qquad \deg(\mathbf{d}) > 2|E| - 2|V| \quad \Rightarrow \quad f(\mathbf{d}) = \deg(\mathbf{d}) + 1 - g,$$

i.e., $f(\mathbf{d})$ eventually equals $\deg(\mathbf{d}) + 1 - g$. Hence $f$ is a Riemann function with offset $C = 1 - g$.

The Baker-Norine formula is an analog of the classical Riemann-Roch formula for algebraic curves or Riemann surfaces; we briefly discuss this in Subsection 2.5.

2.4. **Generalizations of the Baker-Norine Rank.** Many variants of the Baker-Norine rank have been studied. We remark that in literature that generalizes that Baker-Norine rank, e.g., [AM10], one typically studies the function $r = f - 1$ where $f$ is as in (3) for various $\mathcal{N}$, and hence $r$ is initially $-1$ instead of initially $0$.

**Example 2.4.** Amini and Manjunath [AM10] generalized Definition 2.3 by taking $L \subset \mathbb{Z}^n_{\deg 0}$ be any lattice of full rank in $\mathbb{Z}^n_{\deg 0}$ (i.e., rank $n-1$); it this case the definitions of "equivalent," "effective," and of $\mathcal{N}$ in Definition 2.3 carry over; they show that $f$ as in (3) is a Riemann funtion with offset is $1 - g_{\max}(L)$, with $g_{\max}(L)$ as defined on page 5 there. They also give conditions on $L$ so that a Riemann-Roch analog (5) holds; one of their conditions is that all maximal points of $\mathcal{N}$ have the same degree (i.e., $g_{\min} = g_{\max}$ as in [AM10]); they give a second, more technical condition.

To generalize the above examples, let us give some conditions on a subset $\mathcal{N} \subset \mathbb{Z}^n$ which ensure that $f$ in (3) gives a Riemann function.

**Proposition 2.5.** *Let $n \in \mathbb{N}$ and $\mathcal{N} \subset \mathbb{Z}^n$ such that*

*(1) for some $m, m' \in \mathbb{Z}$ we have*

(7) $$\mathbb{Z}^n_{\deg \le m} \subset \mathcal{N} \subset \mathbb{Z}^n_{\deg \le m'},$$

*and*

*(2) setting $M$ to be the largest degree of an element of $\mathcal{N}$, then there exists a $C$ such that if $\mathbf{d} \in \mathbb{Z}^n_{\deg M}$, then then some $\mathbf{d}' \in \mathcal{N} \cap \mathbb{Z}^n_{\deg M}$ has $\|\mathbf{d} - \mathbf{d}'\|_1 \le C$.*

*Then $f$ as in (3) is a Riemann function with offset $-M$.*

*Proof.* Since $\mathbf{d} \in \mathcal{N}$ for $\deg(\mathbf{d}) \le m$, we have that $f$ is initially zero. By induction on $\deg(\mathbf{d})$, we easily show that for any $\mathbf{d}$ with $\deg(\mathbf{d}) > M$, the $L^1$ distance from $\mathbf{d}$ to $\mathbb{Z}_{\le M}$ is at least $\deg(\mathbf{d}) - M$. Hence

(8) $$f(\mathbf{d}) \ge \deg(\mathbf{d}) - M;$$

let us show that equality holds for $\deg(\mathbf{d}) \ge M + Cn$. Say that $\mathbf{d} \in \mathbb{Z}^n$ satisfies $\deg(\mathbf{d}) \ge M + Cn$. Then setting $b = \deg(\mathbf{d}) - M - Cn \ge 0$ we have

$$\widetilde{\mathbf{d}} = \mathbf{d} - C\mathbf{1} - b\mathbf{e}_1$$

has degree $M$; hence for some $\mathbf{d}' \in \mathcal{N} \cap \mathbb{Z}^n_M$ we have

$$\widetilde{\mathbf{d}} - \mathbf{d}' = \mathbf{a}$$

where

$$|a_1| + \cdots + |a_n| \le C;$$

hence $|a_i| \le C$ for all $i$. It follows that setting $\mathbf{a}'$ to be

$$\mathbf{a}' = \mathbf{d} - \mathbf{d}' = \mathbf{d} - (\mathbf{a} + \widetilde{\mathbf{d}}) = C\mathbf{1} + b\mathbf{e}_1 - \mathbf{a},$$

we have $a'_1 = C + a_1 + b$ and for $i \ge 2$, $a'_i = C + a_i$, and hence all $a'_i \ge 0$. Hence the $L^1$ distance of $\mathbf{d}$ to $\mathbf{d}'$ is at most

$$a'_1 + \cdots + a'_n = \deg(\mathbf{d}) - \deg(\mathbf{d}') = \deg(\mathbf{d}) - M,$$

and hence $f(\mathbf{d}) \le \deg(\mathbf{d}) - M$. Hence, (8) holds with equality whenever $\deg(\mathbf{d}) \ge M + Cn$. □

Let us make some further remarks on examples provided by Proposition 2.5.

**Remark 2.6.** Condition (2) of Proposition 2.5 on $\mathcal{N}$ above follows from the following stronger condition: for any $\mathcal{N} \subset \mathbb{Z}^n$, say that $\mathbf{d} \in \mathbb{Z}^n$ is an *invariant translation of* $\mathcal{N}$ if for all $\mathbf{d}' \in \mathbb{Z}^n$, $\mathbf{d}' \in \mathcal{N}$ iff $\mathbf{d}+\mathbf{d}' \in \mathcal{N}$. We easily see that the set, $T = T(\mathcal{N})$ of all invariant translations is a subgroup of the additive group $\mathbb{Z}^n$, and that (7) implies that $T \subset \mathbb{Z}^n_{\deg 0}$. If $T$ is a full rank subgroup of $\mathbb{Z}^n_{\deg 0}$ (i.e., of rank $n-1$), then condition (2) of Proposition 2.5 is automatically satisfied.

**Remark 2.7.** In typical examples $\mathcal{N}$ above is a *downset*, i.e., $\mathbf{d} \in \mathcal{N}$ and $\mathbf{d}' \le \mathbf{d}$ implies that $\mathbf{d}' \in \mathcal{N}$. In this case if the closest point in $\mathcal{N}$ to some $\mathbf{d} \in \mathbb{Z}^n$ is $\mathbf{d}' \in \mathcal{N}$, then clearly (1) $\mathbf{d}' \le \mathbf{d}$, and (2) with $f$ as in (3), $f(\mathbf{d}) = \deg(\mathbf{d} - \mathbf{d}')$; we easily verify the converse, i.e.,

$$f(\mathbf{d}) = \min\{\deg(\mathbf{d} - \mathbf{d}') \mid \mathbf{d}' \in \mathcal{N},\ \mathbf{d}' \le \mathbf{d}\}$$
$$= \min\{\deg(\mathbf{d} - \mathbf{d}') \mid f(\mathbf{d}') = 0\}.$$

Furthermore, if $\mathcal{N}$ is a downset, then for any $i \in [n]$, any path from a $\mathbf{d} \in \mathbb{Z}^n$ to a $\mathbf{d}' \in \mathcal{N}$ translates to a path of the same length from $\mathbf{d} - \mathbf{e}_i$ to $\mathbf{d}' - \mathbf{e}_i$, which again lies in $\mathcal{N}$. Hence if $\mathcal{N}$ is a downset, then $f = f(\mathbf{d})$ as in (3) is a non-decreasing function of $\mathbf{d}$.

**Remark 2.8.** We remark that if $L \subset \mathbb{Z}^n_{\deg 0}$ is not of full rank in Example 2.4, then condition (2) of Proposition 2.5 fails to hold, and we easily see that $f$ in (3) fails to be a Riemann function.

2.5. **Examples Based on Riemann's Theorem.** All the above discussion is based on the classical *Riemann's theorem* and *Riemann-Roch theorem*. However, we use these examples only for illustration, and they are not essential to our discussion of the Baker-Norine rank functions and of most of the rest of this article.

Let $X$ be an algebraic curve over an algebraically closed field $k$, and $K$ be its function field; one understands either (1) $K$ is a finite extension of $k(x)$ where $x$ is an indeterminate (i.e., transcendental) and $X$ is its set of discrete valuations (e.g., [Lan82], Section 1.2), or (2) $X$ is projective curve in the usual sense (e.g., [Har77], Section 4.1), and $K$ is its function field. (For $k = \mathbb{C}$ one can also view $X$ as a compact Riemann surface, and $K$ as its field of meromorphic functions.) To each $f \in K \setminus \{0\}$ one associates the divisor (i.e., Weil divisor) equal to $(f) = \sum_{v \in X} \mathrm{ord}_v(f) v$ [Lan82][1]. For each divisor $D$ one sets

$$L(D) = \{0\} \cup \{f \in K \mid (f) \ge -D\},$$

where we regard $0 \in K$ as having divisor $(0) \ge -D$ for all $D$; this makes $L(D) \subset K$ a $k$-linear subspace, and we set

$$l(D) = \dim_k L(D).$$

For a divisor $D$, we use $\deg(D)$ to denote the sum of the $\mathbb{Z}$-coefficients in $D$. For $f \in K \setminus \{0\}$, $f$ has the same number of zeroes and poles, counted with multiplicity, i.e., $\deg((f)) = 0$. It follows that $l(D) = 0$ when $\deg(D) < 0$. *Riemann's theorem* says that for the *genus* $g \in \mathbb{Z}_{\geq 0}$ of $X$, for any divisor $D$ with $\deg(D)$ sufficiently large,
$$l(D) = \deg(D) + 1 - g.$$
Hence for any points $P_1, \ldots, P_n \in X$ we have

(9) $$f(\mathbf{d}) \stackrel{\text{def}}{=} l(d_1 P_1 + \cdots + d_n P_n)$$

is a Riemann function. The Riemann-Roch formula states that
$$l(D) = l(\omega - D) + \deg(D) + 1 - g$$
where $\omega$ is the *canonical divisor*, i.e., the divisor associated to any 1-form.

---
[1]Here $\mathrm{ord}_v(f)$ is (1) 0 if $f(v)$ is finite and non-zero, (2) the multiplicity of the zero at $v$ if $f(v) = 0$, and (3) minus the multiplicity of the pole at $v$ if $f(v) = \infty$.

**Example 2.9.** Let $K$ be an elliptic curve, i.e., a curve of genus $g = 0$, and $P_1, P_2$ two points of the curve. The Riemann-Roch theorem implies that $f(\mathbf{d}) = 0$ if $\deg(\mathbf{d}) < 0$ and $f(\mathbf{d}) = \deg(\mathbf{d}) - 1$ if $\deg(\mathbf{d}) > 0$. Hence it remains to determine $f(\mathbf{d})$ for $\mathbf{d} = (d_1, -d_1)$ of degree 0, and $f(d_1, -d_1)$ is either 0 or 1. If $P_1 - P_2$ has infinite order in the group law (which, for fixed $P_1$, holds for all but countably many $P_2$), then $f(d_1, -d_1) = 1$ iff $d_1 = 0$; by contrast, if $P_1 - P_2$ has order $r \in \mathbb{N}$, then $f(d_1, -d_1) = 1$ iff $d_1$ is divisible by $r$.

### 2.6. Riemann Functions from other Riemann Functions.

**Example 2.10.** If for some $k, n \in \mathbb{N}$, $f_1, \ldots, f_{2k+1}$ are Riemann functions, then so is
$$f_1 - f_2 + f_3 - \cdots - f_{2k} + f_{2k+1}.$$

One can restrict any Riemann function to a subset of its variables, the others taking fixed values, to get a Riemann function on fewer variables. In [FF] the restriction to two variables is the most important. Let us define the appropriate notation.

**Example 2.11.** Let $f \colon \mathbb{Z}^n \to \mathbb{Z}$ be any Riemann function with $f(\mathbf{d}) = \deg(\mathbf{d}) + C$ for $\deg(\mathbf{d})$ sufficiently large. Then for any distinct $i, j \in [n]$ and $\mathbf{d} \in \mathbb{Z}^n$, the function $f_{i,j,\mathbf{d}} \colon \mathbb{Z}^2 \to \mathbb{Z}$ given as

(10) $$f_{i,j,\mathbf{d}}(a_i, a_j) = f(\mathbf{d} + a_i \mathbf{e}_i + a_j \mathbf{e}_j)$$

is a Riemann function $\mathbb{Z}^2 \to \mathbb{Z}$, and for $a_i + a_j$ large we have

(11) $$f_{i,j,\mathbf{d}}(a_i, a_j) = a_i + a_j + C', \quad \text{where} \quad C' = \deg(\mathbf{d}) + C.$$

We call $f_{i,j,\mathbf{d}}$ a *two-variable restriction* of $f$; we may similarly restrict $f$ to one variable or three or more variables, and any such restriction is clearly a Riemann function.

[It turns out that in [FF], it is important that that $C'$ depends only on $\mathbf{d}$ and not on $i, j$.]

### 2.7. Typical Properties of Riemann Functions.
Let us describe some typical properties of Riemann functions above.

**Definition 2.12.** We say that a function $f \colon \mathbb{Z}^n \to \mathbb{Z}$ is

(1) *slowly growing* if for all $\mathbf{d} \in \mathbb{Z}^n$ and $i \in [n]$ we have
$$f(\mathbf{d}) \leq f(\mathbf{d} + \mathbf{e}_i) \leq f(\mathbf{d}) + 1,$$
and

(2) *p-periodic* for a $p \in \mathbb{N}$ if for all $i, j \in [n]$ and all $\mathbf{d} \in \mathbb{Z}^n$ we have
$$f(\mathbf{d} + p\,\mathbf{e}_i - p\,\mathbf{e}_j) = f(\mathbf{d}).$$

We easily see:
  (1) $f$ in (9) is always slowly growing, but not generally periodic;
  (2) $f$ in (3), then (3) is slowly growing whenever $\mathcal{N}$ is a *downset* (as remarked above);
  (3) in Example 2.4, $f$ is $p$-periodic for any $p$ such that each element of $\mathbb{Z}^n_{\deg 0}/L$ has order divisible by $p$ (hence this holds for $p = |\mathbb{Z}^n_{\deg 0}/L|$);
  (4) in Example 2.11, if $f \colon \mathbb{Z}^n \to \mathbb{Z}$ is either slowly growing or $p$-periodic for some $p$, then the same holds of any restriction of $f$ to two (or any number) of its variables.

## 3. The Weight of a Riemann Function, and Generalized Riemann Functions

In this section we define the *weights* of a Riemann function, a notion central to this article.

Since a Riemann function $\mathbb{Z}^2 \to \mathbb{Z}$ eventually equals $d_1 + d_2 + C$, one may consider that one possible generalization of this notion for a function $\mathbb{Z}^3 \to \mathbb{Z}$ might be a function that eventually equals a polynomial of degree two in $d_1, d_2, d_3$. In fact, most everything we say about Riemann functions hold for a much larger class of functions $\mathbb{Z}^n \to \mathbb{Z}$ which we call *generalized Riemann functions*; this includes all polynomials of $d_1, \ldots, d_n$ of degree $n - 1$, but many more functions.

### 3.1. Weights and Möbuis Inversion.

If $f \colon \mathbb{Z}^n \to \mathbb{Z}$ is initially zero, then there is a unique initially zero $W \in \mathbb{Z}^n \to \mathbb{Z}$ for which

$$(12) \qquad f(\mathbf{d}) = \sum_{\mathbf{d}' \le \mathbf{d}} W(\mathbf{d}'),$$

since we can determine $W(\mathbf{d})$ inductively on $\deg(\mathbf{d})$ set

$$(13) \qquad W(\mathbf{d}) = f(\mathbf{d}) - \sum_{\mathbf{d}' \le \mathbf{d},\ \mathbf{d}' \ne \mathbf{d}} W(\mathbf{d}').$$

Recall from (2) the notation $\mathbf{e}_I$ for $I \subset [n]$.

**Proposition 3.1.** *Consider the operator $\mathfrak{m}$ on functions $f \colon \mathbb{Z}^n \to \mathbb{Z}$ defined via*

$$(14) \qquad (\mathfrak{m}f)(\mathbf{d}) = \sum_{I \subset [n]} (-1)^{|I|} f(\mathbf{d} - \mathbf{e}_I),$$

*and the operator on functions $W \colon \mathbb{Z}^n \to \mathbb{Z}$ that are initially zero given by*

$$(15) \qquad (\mathfrak{s}W)(\mathbf{d}) = \sum_{\mathbf{d}' \le \mathbf{d}} W(\mathbf{d}'),$$

*Then if $f$ is any initially zero function, and $W$ is given by the equation $f = \mathfrak{s}W$ (i.e., $W$ is defined inductively by (13)), then $W = \mathfrak{m}f$.*

The above can be viewed as the Möbius inversion formula for the partial order $\le$ on $\mathbb{Z}^n$.

*Proof.* We have $f(\mathbf{d}) = 0$ whenever $\deg(\mathbf{d}) \le b$ for some $b$, and then (14) shows that $(\mathfrak{m}f)(\mathbf{d}) = 0$ for $\deg(\mathbf{d}) \le b$ as well. Since there is a unique initially zero $W$ with $\mathfrak{s}W = f$, it suffices to show that $\mathfrak{sm}f = f$. Since $f$ is initially zero, for any $\mathbf{d} \in \mathbb{Z}^n$ write $(\mathfrak{sm}f)(\mathbf{d})$ as

$$(\mathfrak{sm}f)(\mathbf{d}) = \sum_{\mathbf{d}' \le \mathbf{d}} \sum_{I \subset [n]} (-1)^{|I|} f(\mathbf{d} - \mathbf{e}_I)$$

which is a double sum of finitely many terms since $f$ is initially zero; hence we may rearrange terms, set $\mathbf{d}'' = \mathbf{d} - \mathbf{e}_I$ and write this double sum as

$$\sum_{\mathbf{d}'' \le \mathbf{d}} f(\mathbf{d}'') \, a_{\mathbf{d}''}, \quad \text{where} \quad a_{\mathbf{d}''} = \sum_{I \text{ s.t. } \mathbf{d}'' + \mathbf{e}_I \le \mathbf{d}} (-1)^{|I|};$$

to compute $a_{\mathbf{d}''}$, setting $J = \{j \in [n] \mid d''_j < d_j\}$, we have

$$\sum_{I \text{ s.t. } \mathbf{d}'' + \mathbf{e}_I \le \mathbf{d}} (-1)^{|I|} = \sum_{I \subset J} (-1)^{|I|}$$

which equals 1 if $J = \emptyset$ and otherwise equals 0. It follows that $a_{\mathbf{d}} = 1$, and for $\mathbf{d}'' \ne \mathbf{d}$, we have $a_{\mathbf{d}''} = 0$. $\square$

**Definition 3.2.** Throughout this article we reserve the symbols $\mathfrak{m}, \mathfrak{s}$ for their meanings in (12) and (14). If $f, W$ are initially zero functions $\mathbb{Z}^n \to \mathbb{Z}$ with $f = \mathfrak{s}W$, we say that $f$ *counts* $W$ and that $W$ is the *weight* of $f$. A function $h \colon \mathbb{Z}^n \to \mathbb{Z}$ is *modular* if $f \in \ker \mathfrak{m}$ (i.e., $\mathfrak{m}f$ is the zero function). We say that $f \colon \mathbb{Z}^n \to \mathbb{Z}$ is a *generalized Riemann function* if

(1) $f$ is initially zero, and
(2) $f$ eventually equals a modular function, i.e., for some $h \in \ker \mathfrak{m}$ we have $f(\mathbf{d}) = h(\mathbf{d})$ for $\deg(\mathbf{d})$ sufficiently large.

3.2. **Weights of Riemann Functions $\mathbb{Z}^2 \to \mathbb{Z}$.** We will be especially interested in Riemann functions $\mathbb{Z}^2 \to \mathbb{Z}$ and their weights $W = \mathfrak{m}f$. It is useful to notice that for such functions we that that for any fixed $d_1$ and $d_2$ sufficiently large,

$$f(d_1, d_2) - f(d_1 - 1, d_2) = 1,$$

and hence, for fixed $d_1$,

(16) $$\sum_{d_2 = -\infty}^{\infty} W(d_1, d_2) = 1,$$

and similarly, for fixed $d_2$ we have

(17) $$\sum_{d_1 = -\infty}^{\infty} W(d_1, d_2) = 1.$$

Viewing $W$ as a two-dimensional infinite array of numbers indexed in $\mathbb{Z} \times \mathbb{Z}$, one can therefore say that $W \colon \mathbb{Z}^2 \to \mathbb{Z}$ is a Riemann weight iff all its "row sums" (16) and all its "column sums" (17) equal one.

### 3.3. Examples and Classification of Generalized Riemann Functions.

At times it is convenient to write $\mathfrak{m}$ using the "downward shift operators," $\mathfrak{t}_i$ for $i \in [n]$, where $\mathfrak{t}_i$ is the operator on functions $\mathbb{Z}^n \to \mathbb{Z}$ given by

$$(\mathfrak{t}_i f)(\mathbf{d}) = f(\mathbf{d} - \mathbf{e}_i); \tag{18}$$

one easily verifies that the $\mathfrak{t}_i$ commute with one another, and that

$$\mathfrak{m} = (1 - \mathfrak{t}_1) \ldots (1 - \mathfrak{t}_n),$$

(where 1 is the identity operator). In particular, it follows that if $f = f(\mathbf{d})$ is independent of its $i$-th variable, then $(1 - \mathfrak{t}_i)f = 0$, and hence $\mathfrak{m} f = 0$. In particular $\mathfrak{m} f = 0$ if (1) $f$ is a sum of functions, each of which is independent in some variable, and, in particular, (2) if $f$ is a polynomial of degree at most $n - 1$. Hence $\deg(\mathbf{d}) + C$ is a modular function for any $n \geq 1$, and hence a Riemann function is, indeed, a generalized Riemann function.

We now characterize modular functions in two different ways.

**Theorem 3.3.** *A function $h \colon \mathbb{Z}^n \to \mathbb{Z}$ is modular iff it can be written as a sum of functions each of which depends on only $n - 1$ of its $n$ variables.*

We postpone its proof to Section 7. The following description of modular functions will be needed when we discuss what we call *Riemann-Roch formulas*.

**Theorem 3.4.** *If $a \in \mathbb{Z}$, $n \in \mathbb{N}$, and $h$ is any integer-valued function defined on $\mathbf{d} \in \mathbb{Z}^n$ with $a \leq \deg(\mathbf{d}) \leq a + n - 1$, then $h$ has a unique extension to a modular function $\mathbb{Z}^n \to \mathbb{Z}$.*

We also postpone the proof of this theorem to Section 7.

According to this theorem, if $h_1, h_2$ are two modular functions, then $h_1$ and $h_2$ are equal whenever they are eventually equal (i.e., $h_1(\mathbf{d}) = h_2(\mathbf{d})$ for $\deg(\mathbf{d})$ sufficiently large), then $h_1 = h_2$. In particular, if $f \colon \mathbb{Z}^n \to \mathbb{Z}$ is a generalized Riemann function, then the modular function $h$ that is eventually equal to $f$ is uniquely determined.

### 3.4. The Weight of the Baker-Norine Rank and Other Functions Initially Equal to $-1$.

Since the Baker-Norine rank and many similar functions are initially equal to $-1$, we make the following convention.

**Definition 3.5.** *If $r \colon \mathbb{Z}^n \to \mathbb{Z}$ is a function that is initially equal to $-1$, by the* weight *of $r$ we mean the function $\mathfrak{m} r$, which clearly equals $\mathfrak{m} f$ with $f = 1 + r$.*

We also note that in the above definition, for any $i \in [n]$ we have $(1 - \mathfrak{t}_i)r = (1 - \mathfrak{t}_i)f$. Hence, as soon as we apply either all of $\mathfrak{m}$, or merely one of its factors $1 - \mathfrak{t}_i$, there is no difference in working with $r$ or $f$. When computing the weight of Baker-Norine type functions, we often use the more suggestive $r_{\mathrm{BN}}$ rather than $f = 1 + r_{\mathrm{BN}}$.

## 4. Riemann-Roch Formulas and Self-Duality

In this section we express Riemann-Roch formulas more simply in terms of the weight of the Riemann function.

**Definition 4.1.** Let $f\colon \mathbb{Z}^n \to \mathbb{Z}$ be a generalized Riemann function, and $h$ the modular function eventually equal to $f$. For $\mathbf{K} \in \mathbb{Z}^n$, the **K**-*dual of* $f$, denoted $f_{\mathbf{K}}^{\wedge}$, refers to the function $\mathbb{Z}^n \to \mathbb{Z}$ given by

$$f_{\mathbf{K}}^{\wedge}(\mathbf{d}) = f(\mathbf{K} - \mathbf{d}) - h(\mathbf{K} - \mathbf{d}). \tag{19}$$

We equivalently write

$$f(\mathbf{d}) - f_{\mathbf{K}}^{\wedge}(\mathbf{K} - \mathbf{d}) = h(\mathbf{d}) \tag{20}$$

and refer to this equation as a *generalized Riemann-Roch formula*.

In particular, if $f$ is a Riemann function with offset $C$, then $h(\mathbf{d}) = \deg(\mathbf{d}) + C$, and (20) means that

$$f(\mathbf{d}) - f_{\mathbf{K}}^{\wedge}(\mathbf{K} - \mathbf{d}) = \deg(\mathbf{d}) + C. \tag{21}$$

The usual Riemann-Roch formulas—the classical one and the Baker-Norine formula—are cases where $f_{\mathbf{K}}^{\wedge} = f$ equals $f$ for some $f, \mathbf{K}$. Hence the above definition is very loose: it says that for any generalized Riemann function, $f$, and any $\mathbf{K} \in \mathbb{Z}^n$, there is always a "generalized Riemann-Roch formula;" we refer to the special cases where $f = f_{\mathbf{K}}^{\wedge}$ for some $\mathbf{K}$ as *self-duality* in Definition 4.4 below.

In Subsection 1.1 we explained some reasons we work with generalized Riemann-Roch formulas; briefly, these reasons are: (1) requiring self-duality would eliminate many interesting Riemann functions, such as the general ones considered by [AM10], and likely some interesting generalized Riemann functions; and (2) self-duality does not behave well under fixing some of the variables of a Riemann function and considering the resulting restriction.

We now give remarks, a theorem, and examples regarding generalized Riemann-Roch formulas.

**Definition 4.2.** If $W\colon \mathbb{Z}^n \to \mathbb{Z}$ is any function and $\mathbf{L} \in \mathbb{Z}^n$, the **L**-*dual weight of* $W$, denoted $W_{\mathbf{L}}^*$ refers to the function given by

$$W_{\mathbf{L}}^*(\mathbf{d}) = W(\mathbf{L} - \mathbf{d}).$$

It is immediate that $(W_{\mathbf{L}}^*)_{\mathbf{L}}^* = W$.

**Theorem 4.3.** *Let* $f\colon \mathbb{Z}^n \to \mathbb{Z}$ *be a generalized Riemann function, and* $W = \mathfrak{m}f$. *Let* $\mathbf{K} \in \mathbb{Z}^n$ *and let* $\mathbf{L} = \mathbf{K} + \mathbf{1}$.

*(1) we have*

$$\mathfrak{m}(f_{\mathbf{K}}^{\wedge}) = (-1)^n W_{\mathbf{L}}^* = (-1)^n (\mathfrak{m}f)_{\mathbf{L}}^*. \tag{22}$$

*(2) $f_{\mathbf{K}}^{\wedge}$ is a generalized Riemann function, and a Riemann function if $f$ is.*
*(3) $(f_{\mathbf{K}}^{\wedge})_{\mathbf{K}}^{\wedge} = f$.*
*(4) $f_{\mathbf{K}}^{\wedge} = f$ iff $W_{\mathbf{L}}^* = (-1)^n W$.*

*Proof.* Proof of (1): applying $\mathfrak{m}$ to (19) we have

$$(\mathfrak{m}(f_{\mathbf{K}}^{\wedge}))(\mathbf{d}) = \sum_{I \subset [n]} (-1)^{|I|} f_{\mathbf{K}}^{\wedge}(\mathbf{d} - \mathbf{e}_I) \tag{23}$$

which, in view of (19), equals

$$\sum_{I \subset [n]} (-1)^{|I|} \Big( f(\mathbf{K} - \mathbf{d} + \mathbf{e}_I) - h(\mathbf{K} - \mathbf{d} + \mathbf{e}_I) \Big). \tag{24}$$

Substituting $J = [n] \setminus I$, for any $g \colon \mathbb{Z}^n \to \mathbb{Z}$ we can write

$$\sum_{I \subset [n]} (-1)^{|I|} g(\mathbf{K} - \mathbf{d} + \mathbf{e}_I) = \sum_{J \subset [n]} (-1)^{n-|J|} g(\mathbf{K} - \mathbf{d} + \mathbf{1} - \mathbf{e}_J)$$

$$= (-1)^n \sum_{J \subset [n]} (-1)^{|J|} g(\mathbf{K} - \mathbf{d} + \mathbf{1} - \mathbf{e}_J) = (-1)^n (\mathfrak{m} g)(\mathbf{K} - \mathbf{d} + \mathbf{1}) = (-1)^n (\mathfrak{m} g)^*_{\mathbf{L}}(\mathbf{d}).$$

Taking $g = f - h$, and using $\mathfrak{m} f = W$ and $\mathfrak{m} h = 0$, we have (24) equals $(-1)^n W^*_{\mathbf{L}}(\mathbf{d})$, and since this also equals (19) we get (22).

Proof of (2): $f$ is a generalized Riemann function iff $W = \mathfrak{m}$ is of finite support, which is equivalent to $W^*_{\mathbf{L}}$ being of finite support; hence $f$ is a generalized Riemann function iff $f^\wedge_{\mathbf{K}}$ is. Moreover, $f$ is a Riemann function iff in addition (20) has $h(\mathbf{d}) = \deg(\mathbf{d}) + C$; in this case (21) with $\mathbf{d}$ replaced with $\mathbf{K} - \mathbf{d}$ is equivalent to

$$f(K - \mathbf{d}) - f^\wedge_{\mathbf{K}}(\mathbf{d}) = h(K - \mathbf{d})$$

for all $\mathbf{d}$, which reversing the sign gives

$$f^\wedge_{\mathbf{K}}(\mathbf{d}) - f(\mathbf{K} - \mathbf{d}) = -h(\mathbf{K} - \mathbf{d}) = -\deg(\mathbf{K} - \mathbf{d}) + C = \deg(\mathbf{d}) + C',$$

where $C' = C - \deg(\mathbf{K})$.

Proof of (3): we may write (22) as

$$f^\wedge_{\mathbf{K}} = \mathfrak{s}(-1)^n (\mathfrak{m} f)^*_{\mathbf{L}},$$

and hence

$$(f^\wedge_{\mathbf{K}})^\wedge_{\mathbf{K}} = \mathfrak{s}(-1)^n (\mathfrak{m} f^\wedge_{\mathbf{K}})^*_{\mathbf{L}} = \mathfrak{s}(-1)^n \bigl((-1)^n W^*_{\mathbf{L}}\bigr)^*_{\mathbf{L}} = \mathfrak{s} W = f.$$

Proof of (4): $f^\wedge_{\mathbf{K}} = f$ (since both functions are initially zero) iff $\mathfrak{m} f^\wedge_{\mathbf{K}} = \mathfrak{m} f$, and by (22) this is equivalent to $(-1)^n W^*_{\mathbf{L}} = W$. $\square$

**Definition 4.4.** We say that a generalized Riemann function $f \colon \mathbb{Z}^n \to \mathbb{Z}$ is *self-dual* if either of the equivalent conditions holds:

(1) for some $\mathbf{K} \in \mathbb{Z}^n$, $f^\wedge_{\mathbf{K}} = f$;
(2) for some $\mathbf{L} \in \mathbb{Z}^n$, $W^*_{\mathbf{L}} = (-1)^n W$.

Let us remark on the uniqueness of $\mathbf{K}$ and $\mathbf{L}$ in the above definition: if $W^*_{\mathbf{L}_1} = W^*_{\mathbf{L}_2}$, it follows that for all $\mathbf{d} \in \mathbb{Z}^n$,

$$W(\mathbf{d}) = \bigl((W^*_{\mathbf{L}_2})^*_{\mathbf{L}_2}\bigr)(\mathbf{d}) = \bigl((W^*_{\mathbf{L}_1})^*_{\mathbf{L}_2}\bigr)(\mathbf{d}) = W^*_{\mathbf{L}_1}(\mathbf{L}_2 - \mathbf{d}) = W(\mathbf{L}_1 - \mathbf{L}_2 + \mathbf{d}),$$

and therefore $W$ is translation invariant by $\mathbf{L}_1 - \mathbf{L}_2$; since $f = \mathfrak{s} W$, and $\mathfrak{s}$ commutes with translation, $f$ is also translation invariant by $\mathbf{L}_1 - \mathbf{L}_2$. Similarly, if $f^\wedge_{\mathbf{K}_1} = f^\wedge_{\mathbf{K}_2}$, then $W^*_{\mathbf{L}_1} = W^*_{\mathbf{L}_2}$ where $\mathbf{L}_j = \mathbf{K}_j + \mathbf{1}$, and $\mathbf{L}_1 - \mathbf{L}_2 = \mathbf{K}_1 - \mathbf{K}_2$, and hence $f$ and $W$ are both translation invariant by $\mathbf{K}_1 - \mathbf{K}_2$. Hence $f$ and $W$ have the same set of invariant translations, $T \subset \mathbb{Z}^n_{\deg 0}$. Hence $\mathbf{K}$ and $\mathbf{L}$ in Definition 4.4 are unique up to a translation by the set $T$.

We remark that the condition $(-1)^n W^*_{\mathbf{L}} = W$ seems to have more direct symmetry than the equivalent condition $f^\wedge_{\mathbf{K}} = f$; furthermore, in the examples of the $W$ that we compute in Sections 5 and 6, the $W$ are very sparse (i.e., mostly 0), and so verifying $(-1)^n W^*_{\mathbf{L}} = W$ seems simpler.

Of course, the classical or Graph Riemann-Roch formulas, in terms of our Definition 4.4, are assertions that self-duality holds in these cases.

**Example 4.5.** The Baker-Norine [BN07] Graph Riemann-Roch theorem for a graph, $G = (V, E)$, with $V = \{v_1, \ldots, v_n\}$ can be stated as
$$r_{\mathrm{BN},G}(\mathbf{d}) - r_{\mathrm{BN},G}(\mathbf{K} - \mathbf{d}) = \deg(\mathbf{d}) + 1 - g,$$
where $g = |E| - |V| + 1$ and $\mathbf{K} = \sum_i \mathbf{e}_i (\deg_G(v_i) - 2)$. Since $f = r_{\mathrm{BN},G} + 1$ is the associated Riemann function, the left-hand-side above also equals $f(\mathbf{d}) - f_K^{\wedge}(\mathbf{K} - \mathbf{d})$, and hence $f = f_K^{\wedge}$ is self-dual.

**Example 4.6.** Amini and Manjunath [AM10] give conditions for $f$ as in (3) with $\mathcal{N}$ as in Example 2.4 to satisfy self-duality. The first is that all maximal points of $\mathcal{N}$ have the same degree ($g_{\min} = g_{\max}$ in [AM10]); the second is more technical. However, to us these Riemann functions seem interesting to study whether or not self-duality holds.

## 5. The Weight of Two Vertex Graphs and Riemann Functions of Two Variables

In this section we prove the following theorem.

**Theorem 5.1.** *Let $G$ be a graph on two vertices, $v_1, v_2$ with $r \geq 1$ edges joining $v_1$ and $v_2$. Let $r_{\mathrm{BN}} \colon \mathbb{Z}^2 \to \mathbb{Z}$ be the Baker-Norine rank, let $f = 1 + r_{\mathrm{BN}}$, i.e., $f$ is as in (3) in Definition 2.3. Then $\mathbf{d}$ is in the image of the Laplacian iff $\mathbf{d}$ is an integral multiple of $(r, -r)$. Let $W = \mathfrak{m} f$ be the weight of $f$. Then*
$$W(0, 0) = W(1, 1) = \ldots = W(r-1, r-1) = 1;$$
*furthermore $W(\mathbf{d}) = 1$ if $\mathbf{d}$ is equivalent to one of $(i, i)$ with $i = 0, \ldots, r-1$, and otherwise $W(\mathbf{d}) = 0$.*

### 5.1. Perfect Matchings and Slowly Growing Riemann Functions.
In this subsection we make some remarks on weights that we call "perfect matchings."

**Definition 5.2.** Let $W$ be a function $\mathbb{Z}^2 \to \mathbb{Z}$ that is initially and eventually zero. We say that $W$ is a *perfect matching* if there exists a permutation (i.e., a bijection) $\pi \colon \mathbb{Z} \to \mathbb{Z}$ such that

$$(25) \qquad W(i, j) = \begin{cases} 1 & \text{if } j = \pi(i), \text{ and} \\ 0 & \text{otherwise.} \end{cases}$$

It follows that for $\pi$ as above, $\pi(i) + i$ is bounded above and below, since $W$ is initially and eventually 0. Of course, if $W$ is $r$-periodic, i.e., for all $\mathbf{d} \in \mathbb{Z}^2$, $W(\mathbf{d}) = W(\mathbf{d} + (r, -r))$, then $\pi$ is *skew-periodic* in the sense that $\pi(i+r) = \pi(i) - r$ for all $i \in \mathbb{Z}$.

**Proposition 5.3.** *Let $f \colon \mathbb{Z}^2 \to \mathbb{Z}$ be a slowly growing Riemann function, i.e., for $i = 1, 2$ and any $\mathbf{d} \in \mathbb{Z}^2$ we have*
$$f(\mathbf{d}) \leq f(\mathbf{d} + \mathbf{e}_i) \leq f(\mathbf{d}) + 1.$$
*Let $W = \mathfrak{m} f$ be the weight of $f$. Then $W$ takes only the values $0$ and $\pm 1$. Furthermore, for any $\mathbf{d} \in \mathbb{Z}^2$, let $a = f(\mathbf{d})$*

$$(26) \quad W(\mathbf{d}) = 1 \iff f(\mathbf{d} - \mathbf{e}_1) = f(\mathbf{d} - \mathbf{e}_2) = f(\mathbf{d} - \mathbf{e}_1 - \mathbf{e}_2) = a - 1,$$

*and*

$$(27) \quad W(\mathbf{d}) = -1 \iff f(\mathbf{d} - \mathbf{e}_1) = f(\mathbf{d} - \mathbf{e}_2) = a = f(\mathbf{d} - \mathbf{e}_1 - \mathbf{e}_2) + 1.$$

*We say that $f$ is* supermodular *when $W(\mathbf{d}) \geq 0$ for all $0$; in this case $W$ is a perfect matching.*

*Proof.* For $\mathbf{d} \in \mathbb{Z}^2$, let $a = f(\mathbf{d})$. Then $f(\mathbf{d} - \mathbf{e}_1 - \mathbf{e}_2)$ is between $a - 2$ and $a$, since $f$ is slowly growing. We proceed by a case analysis:

(1) if $f(\mathbf{d} - \mathbf{e}_1 - \mathbf{e}_2) = a - 2$, then $f(\mathbf{d} - \mathbf{e}_1)$ differs by at most 1 from both $a$ and $a - 2$, and hence $f(\mathbf{d} - \mathbf{e}_1) = a - 1$; similarly $f(\mathbf{d} - \mathbf{e}_2) = a - 1$, and so $W(\mathbf{d}) = 0$.
(2) if $f(\mathbf{d} - \mathbf{e}_1 - \mathbf{e}_2) = a$, then since $f$ is non-decreasing we have $f(\mathbf{d} - \mathbf{e}_i) = a$ for $i = 1, 2$, and hence $W(\mathbf{d}) = 0$;
(3) if $f(\mathbf{d} - \mathbf{e}_1 - \mathbf{e}_2) = a - 1$, then since $f$ is non-decreasing we have that for each $i = 1, 2$, $f(\mathbf{d} - \mathbf{e}_i)$ is either $a$ or $a - 1$; this gives four cases to check, which imply (26) and (27).

If $W$ never takes the value $-1$, then (16) implies that for each $d_1$ there is a unique $d_2$ with $W(d_1, d_2) = 1$, so setting $\pi(d_1) = d_2$ gives a map $\pi \colon \mathbb{Z} \to \mathbb{Z}$; then (17) implies that $\pi$ has an inverse. $\square$

*Proof of Theorem 5.1.* The rows of the Laplacian of $G$ are $(r, -r)$ and $(-r, r)$, and hence the image, $L$, of the Laplacian equals the integer multiples of $(r, -r)$.

First let us prove that $f$ is supermodular by a case analysis: indeed,

(1) if $f(\mathbf{d}) = 0$, then $f(\mathbf{d}') = 0$ for $\mathbf{d}' \leq \mathbf{d}$ and hence $W(\mathbf{d}) = 0$;
(2) if $f(\mathbf{d}) \geq 1$, then there is a path from $\mathbf{d}$ to $\mathcal{N}$ as in (3) of positive length through the points of $\mathbb{Z}^2$, and hence for some $i = 1, 2$ we have $f(\mathbf{d} - \mathbf{e}_i) = f(\mathbf{d}) - 1$; then Proposition 5.3 implies that $W(\mathbf{d}) \geq 0$.

It follows that $W$ is a perfect matching, and hence $W$ is given by (25) for some perfect matching $\pi$; since $f$ is $r$-periodic, it suffices to determine $\pi(i)$ for $i = 0, 1, \ldots, r - 1$. Let us do so by finding some values of $f$.

Since $(0, 0) \in L$, we have $f(0, 0) = 1$, and for all $i \geq 0$, $f(i, 0) \geq 1$. But $(i, 0) - \mathbf{e}_2$ cannot be effective for $i \leq r - 1$, since then for some $m \in \mathbb{Z}$ we would have $(i, -1) \geq m(r, -r)$, which implies both $m \leq i/r < 1$ and $m \geq 1/r > 0$, which is impossible. Hence for $0 \leq i \leq r - 1$ we have $f(i, 0) = 1$.

On the other hand, we can prove that for $i \geq 0$ we have $f(i, i) \geq i + 1$, using induction on $i$: for $i = 0$ we have $f(0, 0) = 1$, and for the inductive claim with $i \geq 1$, since $(i, i)$ is effective we have

$$f(i, i) = 1 + \max\bigl(f(i - 1, i), f(i, i - 1)\bigr) \geq 1 + f(i - 1, i - 1) \geq 1 + i$$

by the inductive hypothesis.

For $0 \leq i \leq r - 1$, since $f(i, 0) = 1$ and $f(i, i) \geq i + 1$, the fact that $f$ is slowly growing implies that $f(i, j) = j + 1$ for $0 \leq j \leq i$. Similarly, for such $i, j$ with $0 \leq i \leq j$, $f(i, j) = i + 1$.

Using this, it follows that for $i = 0, \ldots, r - 1$ we have

$$W(i, i) = f(i, i) - 2f(i, i - 1) + f(i - 1, i - 1) = i - 2(i - 1) + i - 1 = 1.$$

It follows that $\pi(i) = i$ for $0 \leq i \leq r - 1$, and the theorem follows. $\square$

Notice that this computation proves the Riemann-Roch formula in this case: this computation shows that $W = W_\mathbf{L}^*$ for $\mathbf{L} = (r - 1, r - 1)$. Hence $f = f_\mathbf{K}^\wedge$ for $\mathbf{K} = (r - 2, r - 2)$, and therefore

$$f(\mathbf{d}) - f(\mathbf{K} - \mathbf{d}) = \deg(\mathbf{d}) + C$$

for some $C$. Taking $\mathbf{d} = 0$ and using $f(0,0) = 1$ we get
$$1 - f(\mathbf{K}) = C,$$
and taking $\mathbf{d} = \mathbf{K}$ we get
$$f(\mathbf{K}) - 1 = \deg(\mathbf{K}) + C = 2(r-2) + C;$$
adding these last two equations, the $f(\mathbf{K})$ cancels and we get $0 = 2(r-2) + 2C$, and so $C = 2 - r$ is the offset. Hence
$$f(\mathbf{d}) - f(\mathbf{K} - \mathbf{d}) = \deg(\mathbf{d}) - r + 2.$$

## 6. The Weight of the Riemann-Roch Rank of the Complete Graph and Related Graphs

The point of this subsection is to give a self-contained computation of the remarkably simple and sparse weight function of the Baker-Norine rank for the complete graph.

Our proof uses many standard ideas in the graph Riemann-Roch literature [BN07, Bac17, AM10, CB13], but also one rather ingenious idea of Cori and Le Borgne [CB13].

### 6.1. Proof Overview and Computer-Aided Computations.
Our analysis of the weights for the complete graph and the resulting formula of the Baker-Norine function is based on seeing some remarkable patterns in computer-aided computation. Explaining this also serves as an overview for our proofs below, and motivates the notation that we introduce.

Let $G$ be a graph on $n$-vertices ordered $v_1, \ldots, v_n$. To compute the Baker-Norine function, $r_{\mathrm{BN}}$ of a graph (and the resulting weight, $W$), we note tht $r_{\mathrm{BN}}(\mathbf{d}) = -1$ if $\deg(\mathbf{d}) < 0$; it suffices to compute $r_{\mathrm{BN}}(\mathbf{d})$ on $\mathbb{Z}^n_{\deg 0}$, then on $\mathbb{Z}^n_{\deg 1}$, then $\mathbb{Z}^n_{\deg 2}$, etc. Since $r_{\mathrm{BN}}$ and $W$ are invariant under the image of the Laplacian, $\Delta_G$, it suffices to determine the value of $r_{\mathrm{BN}}$ on a set of representatives of
$$\mathrm{Pic}_i(G) = \mathbb{Z}^n_{\deg i}/\mathrm{Image}(\Delta_G)$$
for $i = 0, 1, \ldots$. To do so, it is natural to: find a set of "convenient coordinates" for $\mathrm{Pic}_0(G) = \mathbb{Z}^n_{\deg 0}/\mathrm{Image}(\Delta_G)$, meaning a set $\mathcal{B}$ and a bijection $\iota \colon \mathcal{B} \to \mathrm{Pic}_0(G)$ such that the computations below are easy to do for $i = 0, 1, \ldots$, namely:

(1) for all $\mathbf{b} \in \mathcal{B}$, determine if $\iota(\mathbf{b}) + i\mathbf{e}_n$ is not effective, i.e., if $r_{\mathrm{BN}}(\iota(\mathbf{b}) + i\mathbf{e}_n) = -1$; and
(2) for all other $\mathbf{b} \in \mathcal{B}$ we compute $r_{\mathrm{BN}}(\mathbf{b} + i\mathbf{e}_n)$ via the formula
$$r_{\mathrm{BN}}(b + i\mathbf{e}_n) = 1 + \min_{j \in [n]} r_{\mathrm{BN}}(\mathbf{b} + i\mathbf{e}_n - \mathbf{e}_j);$$
hence we need a reasonably fast algorithm to determine the element of $\mathcal{B}$ that is equivalent to $\iota^{-1}(\mathbf{b} + \mathbf{e}_n - \mathbf{e}_j)$. [We are finished when $i \geq \deg(\mathbf{L})$ where $\mathbf{L} = \mathbf{K} + \mathbf{1}$ where $K$ is the Baker-Norine canonical divisor, and hence when $i \geq 2(|E| - |V|) + |V| = 2|E| - |V|$; we may use $W = (-1)^n W^*_{\mathbf{L}}$ to finish when $i \geq |E| + (1 - |V|)/2$.]

Of course, one can replace $\mathbf{e}_n$ above by any of $\mathbf{e}_1, \ldots, \mathbf{e}_{n-1}$, or, more generally, any element of $\mathbb{Z}^n$ of degree 1; our choice of $\mathbf{e}_n$ is convenient for the representatives of $\mathcal{B}$ below.

It turns out that there is a very convenient choice for $\mathcal{B}$ suggested in [CB13]: namely, we give their proof that every element of $\mathbb{Z}^n$ is equivalent to a unique element of $\mathcal{A}$ given by

$$\mathcal{A} = \{\mathbf{a} \mid a_1, \ldots, a_{n-2} \in \{0, \ldots, n-1\}, a_{n-1} = 0\},$$

i.e., some element of the form

$$(a_1, \ldots, a_n) \in \mathcal{A} = \{0, \ldots, n-1\}^{n-2} \times \{0\} \times \mathbb{Z} \subset \mathbb{Z}^n$$

The only problem is that the group law in $\mathrm{Pic}(K_n)$ is a bit tricky to write down, since if $\mathbf{a}, \mathbf{a}' \in \mathcal{A}$, then the element of $\mathcal{A}$ that is equivalent to $\mathbf{a} + \mathbf{a}'$ has, for all $i \le n-2$, its $i$-th coordinate equal to $(a_i + a_i') \bmod n$, but the $n$-th coordinate needs to take into account the number of $i$ such that $a_i + a_i' \ge n$. In other words, the addition law on the first $n-2$ coordinates of $\mathcal{A}$ is that of $(\mathbb{Z}/n\mathbb{Z})^{n-2}$ (and the $(n-1)$-th coordinate is always 0), but addition on the $n$-th coordinate depends on the first $n-2$ coordinates; in other words, the addition law on $\mathcal{A}$ induced by the law on Pic gives an isomorphism between $\mathcal{A}$ and a semidirect product $(\mathbb{Z}/n\mathbb{Z})^{n-2} \ltimes \mathbb{Z}$.

Of course, since $\mathcal{A} \subset \mathbb{Z}^n$, this type of complicated addition law cannot be helped: the order of any nonzero element of $\mathbb{Z}^n$ is infinite, whereas the order of each element in $\mathrm{Pic}_0$ is finite; hence if $\mathrm{Pic}_0$ is nontrivial (or, equivalently, $G$ is not a tree), then no set of representatives of Pic can have a simple addition law.

To get a simpler addition law, we define a second set of coordinates: namely, we set $\mathcal{B} = \{0, \ldots, n-1\}^{n-2}$, we define $\iota \colon \mathcal{B} \to \mathrm{Pic}_0$ via

$$\iota \mathbf{b} = (b_1, \ldots, b_{n-2}, 0, -b_1 - \cdots - b_{n-2}) \in \mathbb{Z}^n_{\deg 0}.$$

In order to avoid writing $\iota$ all the time, for $(\mathbf{b}, i) \in \mathcal{B} \times \mathbb{Z}$ we set

$$\langle \mathbf{b}, i \rangle = \iota(\mathbf{b}) + i\mathbf{e}_n,$$

which equals

$$(b_1, \ldots, b_{n-2}, 0, i - b_1 - \cdots - b_{n-2}) \in \mathbb{Z}^n_{\deg i}.$$

Hence we leave the first $n-1$ coordinates as is in $\mathcal{A}$, but we form $\langle \mathbf{b}, i \rangle$ to have degree $i$. In this way

$$\langle \mathbf{b}, i \rangle + \langle \mathbf{b}', i' \rangle$$

has degree $i+i'$, has $(n-1)$-th coordinate 0, and has the first $n-2$ coordinates given by addition in $(\mathbb{Z}/n\mathbb{Z})^{n-2}$; hence the addition law in Pic in the second coordinates $(\mathbf{b}, i)$, is just addition on $(\mathbb{Z}/n\mathbb{Z})^{n-2} \times \mathbb{Z}$.

The theorems we give below simply reflect the patterns that we saw, namely: we first noticed that the weights $W = \mathfrak{m} r_{\mathrm{BN}}$ for the complete graph were very sparse, i.e., mostly 0's, and the non-zero values of $W$ followed a simple pattern. Then, since

$$\mathfrak{m} = (1 - \mathfrak{t}_1) \ldots (1 - \mathfrak{t}_n)$$

(recall that $\mathfrak{t}_i$ is the "downward shift operator" given in (18)), we tried computing some subset of the $1 - \mathfrak{t}_i$ applied to $r_{\mathrm{BN}}$ to find a simple pattern. After a number of unsuccessful attempts, we discovered that $(1 - \mathfrak{t}_{n-1}) r_{\mathrm{BN}}$ had a remarkably simple pattern, namely that for small $n$,

$$(1 - \mathfrak{t}_{n-1}) r_{\mathrm{BN}}(\langle \mathbf{b}, i \rangle) = \begin{cases} 1 & \text{if } b_1 + \cdots + b_n \le i \\ 0 & \text{otherwise.} \end{cases}$$

From this one also easily sees the pattern

$$(1 - \mathfrak{t}_n)(1 - \mathfrak{t}_{n-1})r_{\mathrm{BN}}(\langle \mathbf{b}, i \rangle) = \begin{cases} 1 & \text{if } b_1 + \cdots + b_n = i \\ 0 & \text{otherwise.} \end{cases}$$

The rest of this section is devoted to proving that these patterns above, which we observed for small $n$, indeed hold for all $n$. Our starting point for the proof requires some important techniques of [CB13], which are more simply stated in terms of the representatives $\mathcal{A}$ of $\mathrm{Pic}(K_n) = \mathbb{Z}^n/\mathrm{Image}(\Delta_{K_n})$ used by used in [CB13].

6.2. **Maximal Decrease.** The following is a standard tool used in studying the graph Riemann-Roch rank, used by Baker-Norine [BN07] and many subsequent papers. It is valid in the general setting of (3) when $\mathcal{N}$ is a downset.

Recall from Definition 2.12 that $f \colon \mathbb{Z}^n \to \mathbb{Z}$ if for all $j \in [n]$ and $\mathbf{d} \in \mathbb{Z}^n$ we have

$$f(\mathbf{d}) \le f(\mathbf{d} + \mathbf{e}_j) \le f(\mathbf{d}) + 1.$$

If so, an easy induction argument (on $\deg(\mathbf{d} - \mathbf{d}')$) shows that if $\mathbf{d}', \mathbf{d} \in \mathbb{Z}^n$ with $\mathbf{d}' \le \mathbf{d}$, then

(28) $$f(\mathbf{d}') \ge f(\mathbf{d}) - \deg(\mathbf{d} - \mathbf{d}').$$

**Definition 6.1.** Let $f \colon \mathbb{Z}^n \to \mathbb{Z}$ be slowly growing. Let $\mathbf{d}', \mathbf{d} \in \mathbb{Z}^n$ with $\mathbf{d}' \le \mathbf{d}$. We say that $f$ *is maximally decreasing from* $\mathbf{d}$ *to* $\mathbf{d}'$ if equality holds in (28), or equivalently

$$f(\mathbf{d}) = f(\mathbf{d}') + \deg(\mathbf{d} - \mathbf{d}').$$

The following is Lemma 5 of [CB13], but is used in most papers we have seen involving the Baker-Norine rank, e.g., [BN07, Bac17, AM10].

**Proposition 6.2.** *Let $f \colon \mathbb{Z}^n \to \mathbb{Z}$ be slowly growing. Then for any $\mathbf{d}'', \mathbf{d}', \mathbf{d} \in \mathbb{Z}^n$, $f$ is maximally decreasing from $\mathbf{d}$ to $\mathbf{d}''$ iff it is maximally decreasing from both $\mathbf{d}$ to $\mathbf{d}'$ and from $\mathbf{d}'$ to $\mathbf{d}''$.*

The proof is immediate from the fact that the two inequalities

$$f(\mathbf{d}) - f(\mathbf{d}') \le \deg(\mathbf{d} - \mathbf{d}'),$$
$$f(\mathbf{d}') - f(\mathbf{d}'') \le \deg(\mathbf{d}' - \mathbf{d}'')$$

both hold with equality iff their sum does, and their sum is

$$f(\mathbf{d}) - f(\mathbf{d}'') \le \deg(\mathbf{d} - \mathbf{d}').$$

We remark that $f$ is slowly growing whenever it is of the form (3) where $\mathcal{N}$ is a downset such that $\mathbb{Z}^n_{\deg \le m} \subset \mathcal{N}$ for some $m$ (so that $f$ takes on finite values). We also remark that in this case $\mathbf{d} \in \mathbb{Z}^n$, and $\mathbf{d}'' \in \mathcal{N}$ is such that

$$\|\mathbf{d} - \mathbf{d}''\| = \min_{\mathbf{d}' \in \mathcal{N}} \|\mathbf{d} - \mathbf{d}'\|,$$

then $f$ is maximally decreasing from $\mathbf{d}$ to $\mathbf{d}''$.

### 6.3. A Generalization of a Fundamental Lemma of Cori and Le Borgne.
Next we give an elegant and rather ingenious observation of [CB13] (half of the proof of Proposition 10 there) that is the starting point of their (and our) study the Baker-Norine rank for the complete graph; we state their observation in slightly more general terms.

**Lemma 6.3.** *Fix $n \in \mathbb{N}$, and let $K_n = (V, E)$ be the complete graph on vertex set $V = [n]$, i.e., $E$ consists of exactly one edge joining any two distinct vertices. Consider the Baker-Norine rank $r_{\mathrm{BN}} \colon \mathbb{Z}^n \to \mathbb{Z}$ on $K_n$. If $\mathbf{a} \geq \mathbf{0}$ then*

$$(29) \qquad a_{n-1} = 0 \quad \Rightarrow \quad r_{\mathrm{BN}}(\mathbf{a}) = r_{\mathrm{BN}}(\mathbf{a} - \mathbf{e}_{n-1}) + 1.$$

Of course, by symmetry (29) holds with both occurrences of $n - 1$ replaced by any $j \in [n]$.

*Proof.* Since $\mathbf{a} \geq \mathbf{0}$, $r_{\mathrm{BN}}(\mathbf{a}) \geq 0$, and hence $r_{\mathrm{BN}}$ is maximally decreasing from $\mathbf{a}$ to $\mathbf{a} - \mathbf{b}$ for some $\mathbf{b} \geq \mathbf{0}$ with $r_{\mathrm{BN}}(\mathbf{a} - \mathbf{b}) = -1$. Since $r_{\mathrm{BN}}(\mathbf{a} - \mathbf{b}) = -1$, we must have $a_j - b_j \leq -1$ for some $j \in [n]$; fix any such $j$. Then $b_j \geq a_j + 1 \geq 1$; setting $\mathbf{a}' = \mathbf{a} - b_j \mathbf{e}_j$ we have

$$\mathbf{a} - \mathbf{b} \leq \mathbf{a}' \leq \mathbf{a},$$

and hence $r_{\mathrm{BN}}$ is maximally decreasing from $\mathbf{a}$ to $\mathbf{a}'$. But the vector

$$(30) \qquad \mathbf{a}'' = \mathbf{a} - a_j \mathbf{e}_j - (b_j - a_j) \mathbf{e}_{n-1}$$

is merely the vector $\mathbf{a}'$ followed by an exchange of the $(n-1)$-th and $j$-th coordinates (if $j = n - 1$, then $\mathbf{a}'' = \mathbf{a}'$). Hence $\mathbf{a}'', \mathbf{a}'$ have the same degree and same value of $r_{\mathrm{BN}}$; hence $f$ is also maximally decreasing from $\mathbf{a}$ to $\mathbf{a}''$. Since $b_j - a_j \geq 1$, (30) implies

$$\mathbf{a}'' \leq \mathbf{a} - \mathbf{e}_{n-1} \leq \mathbf{a};$$

since $f$ is maximally decreasing from $\mathbf{a}$ to $\mathbf{a}''$, $f$ is maximally decreasing from $\mathbf{a}$ to $\mathbf{a} - \mathbf{e}_{n-1}$ as well, and hence (29) holds. □

**Remark 6.4.** If $n, m \in \mathbb{N}$, we use $K_n^m = (V, E)$ to denote the graph with $V = [n]$ and $m$ edges between any two vertices (so $K_n^1 = K_n$). Then $r_{\mathrm{BN}, K_n^m}(\mathbf{d})$ is again a symmetric function of its variables $(d_1, \ldots, d_n) = \mathbf{d}$, and the same argument shows that for any $b \in \mathbb{Z}_{\geq 0}$, $\mathbf{a} \geq b\mathbf{1}$ and $a_{n-1} = b$ implies that $f(\mathbf{d}) = f(\mathbf{d} - \mathbf{e}_{n-1}) + 1$. We believe it is possible to use this observation, specifically for $b = m$, to give an analog of Theorem 6.17 below regarding $K_n^m$.

### 6.4. The First Coordinates for Pic, D'après Cori-Le Borgne.
Let us recall some more standard graph Riemann-Roch terminology (see, e.g., [BN07, CB13], and then give our first set of coordinates for the *Picard group* of a graph. These coordinates are those found in the Algorithm at the end of Section 2.1 of [CB13].

Recall $\mathbb{Z}^n_{\deg i}$ consists of the elements of $\mathbb{Z}^n$ of degree $i$. Recall [BN07] the *Picard group* of a graph, $G$, with $n$ vertices $v_1, \ldots, v_n$ is defined as

$$\mathrm{Pic}(G) = \mathbb{Z}^n / \mathrm{Image}(\Delta_G);$$

since $\mathrm{Image}(\Delta_G)$ consists entirely of vectors of degree 0, $\mathrm{Pic}(G)$ is the union over $i \in \mathbb{Z}$ of

$$(31) \qquad \mathrm{Pic}_i(G) = \mathbb{Z}^n_{\deg i} / \mathrm{Image}(\Delta_G).$$

It is known that for all $i$, $|\mathrm{Pic}_i(G)|$ equals $(1/n) \det'(\Delta_G)$, where $\det'$ denotes the product of the nonzero eigenvalues of $\Delta_G$ (and Kirchoff's theorem says that this is the number of unrooted spanning trees of $G$). For $G = K_n$ it is a standard fact that this number of trees is $n^{n-2}$, i.e.,

(32) $$|\operatorname{Pic}_i(K_n)| = n^{n-2}.$$

Next we pick a convenient set of representatives for each class in $\mathbb{Z}^n/\operatorname{Image}(\Delta_{K_n})$.

**Notation 6.5.** For any $n \in \mathbb{N}$, we let

(33) $$\mathcal{A} = \mathcal{A}(n) = \{\mathbf{a} \in \mathbb{Z}^n \mid a_1, \ldots, a_{n-2} \in \{0, \ldots, n-1\}, a_{n-1} = 0\}$$
$$= \{0, \ldots, n-1\}^{n-2} \times \{0\} \times \mathbb{Z}$$

(we usually simply write $\mathcal{A}$ since $n$ will be understood and fixed); in addition, for $i \in \mathbb{Z}$, we use $\mathcal{A}_{\deg i}$ to denote the set

$$\mathcal{A}_{\deg i} \stackrel{\text{def}}{=} \mathcal{A} \cap \mathbb{Z}^n_{\deg i} = \{\mathbf{a} \in \mathcal{A} \mid \deg(\mathbf{a}) = i\}.$$

In the above notation, note that

$$\mathbf{a} \in \mathcal{A}_{\deg i} \iff a_n = i - a_1 - \cdots - a_{n-2}$$

and hence

(34) $$\mathbf{a} \in \mathcal{A}_{\deg i} \Rightarrow \Big(a_n \geq 0 \iff a_1 + \cdots + a_{n-2} \leq i\Big)$$

(35) $$\mathbf{a} \in \mathcal{A}_{\deg i} \Rightarrow \Big(a_n = 0 \iff a_1 + \cdots + a_{n-2} = i\Big)$$

**Lemma 6.6.** *Fix $n \in \mathbb{N}$, and let $K_n = (V, E)$ be the complete graph on vertex set $V = [n]$. Then for all $\mathbf{d} \in \mathbb{Z}^n$ there exists a unique $\mathbf{a} \in \mathcal{A} = \mathcal{A}(n)$ with $\mathbf{d} \sim \mathbf{a}$ (i.e., $\mathbf{d} - \mathbf{a} \in \operatorname{Image}(\Delta_{K_n})$), given by: for $j \in [n-2]$, $a_j = (d_j - d_{n-1}) \bmod n$, i.e., $a_j$ is the element of $\{0, \ldots, n-1\}$ congruent to $d_j - d_{n-1}$ modulo $n$, $a_{n-1} = 0$, and $a_n = \deg(\mathbf{d}) - a_1 - \cdots - a_{n-2}$.*

*Proof.* Existence is shown in "Algorithm" at the end of Section 2.1 of [CB13]: we note that the image of $\Delta_G$ contains $(1, \ldots, 1, 1-n)$ and, for any $j \in [n]$, $n(\mathbf{e}_j - \mathbf{e}_n)$. For any $\mathbf{d}$ we get an equivalent vector with $(n-1)$-th coordinate 0 by subtracting multiples of $(1, \ldots, 1, 1-n)$; then we find an equivalent vector with the first $n-2$ coordinates between 0 and $n-1$ by subtracting multiples of $n(\mathbf{e}_j - \mathbf{e}_n)$ for $j \in [n-2]$.

Note that the above algorithm determines a map $\mu \colon \mathbb{Z}^n \to \mathcal{A}$ that such that

(36) $$\forall \mathbf{d} \in \mathbb{Z}^n, \quad \mathbf{d} \sim \mu(\mathbf{d}),$$

i.e., $\mathbf{d}$ and $\mu(\mathbf{d})$ are equivalent modulo $\operatorname{Image}(K_n)$.

To prove that each $\mathbf{d}$ is equivalent to a unique element of $\mathcal{A}$, we need to show that if $\mathbf{a}, \mathbf{a}' \in \mathcal{A}$ are equivalent, i.e., $\mathbf{a} - \mathbf{a}' \in \operatorname{Image}(\Delta_{K_n})$, then we must have $\mathbf{a} = \mathbf{a}'$. Note that if $\mathbf{a}, \mathbf{a}'$ are equivalent, then they have the same degree and hence both lie in $\mathcal{A}_{\deg i}$ for the same $i$. Hence it suffices to show that each element of $\mathcal{A}_{\deg i}$ is in a distinct class of $\operatorname{Pic}_i(K_n)$. Let us rephrase this condition.

Note that since $\mathcal{A}_{\deg i} \subset \mathbb{Z}^n_{\deg i}$, the quotient map

$$\mathbb{Z}^n_{\deg i} \to \mathbb{Z}^n_{\deg i}/\operatorname{Image}(\Delta_{K_n}) = \operatorname{Pic}_i(K_n)$$

restricts to a map

$$\nu_i \colon \mathcal{A}_{\deg i} \to \operatorname{Pic}_i(K_n).$$

To show that each element of $\mathcal{A}_{\deg i}$ is in its own class of $\operatorname{Pic}_i(K_n)$ simply means that $\nu_i$ is injective. Let us prove this.

So fix an $i \in \mathbb{Z}$. Choosing a set of representatives, $\mathcal{P}_i \subset \mathbb{Z}^n_i$ for $\operatorname{Pic}_i$; in view of (36), $\mu$ restricted to $\mathcal{P}_i$ gives a map of sets $\mu|_{\mathcal{P}_i} \colon \mathcal{P}_i \to \mathcal{A}_{\deg i}$ that takes each element in the domain to a vector equivalent to it; hence this gives a map of sets $\mu_i \colon \mathrm{Pic}_i \to \mathcal{A}_{\deg i}$ such that $\mu_i$ takes each $p \in \mathrm{Pic}_i$ to an element that lies in $p$. It follows that the map $\nu_i \mu_i$ is the identity map on $\mathrm{Pic}_i$.

But we easily see that $\mathcal{A}_{\deg i}$ has size $n^{n-2}$, since if $\mathbf{a} = (a_1, \ldots, a_n) \in \mathcal{A}_{\deg i}$ then $a_1, \ldots, a_{n-2} \in \{0, \ldots, n-1\}$, and any $a_1, \ldots, a_{n-2} \in \{0, \ldots, n-1\}$ determine the values of $a_{n-1}, a_n$, namely

$$a_{n-1} = 0, \quad a_n = i - a_1 - \cdots - a_{n-2}.$$

Since $\nu_i \mu_i$ is the identity map on $\mathrm{Pic}_i$, and this map factors through the set $\mathcal{A}_{\deg i}$ of the same size, both $\nu_i$ and $\mu_i$ must be bijections. Hence $\nu_i$ is an injection, which proves the desired uniqueness property. $\square$

Here is how we often use the above theorem.

**Corollary 6.7.** *Fix an $n \in \mathbb{N}$. For each $i \in \mathbb{Z}$, $\mathcal{A}_{\deg i}$ is a set of representatives of the classes $\mathrm{Pic}_i(K_n)$ in $\mathbb{Z}^n_{\deg i}$. Similarly, for any $\mathbf{d} \in \mathbb{Z}^n$, as $\mathbf{a}$ ranges over $\mathcal{A}_{\deg i}$, $\mathbf{a} - \mathbf{d}$ ranges over a set of representatives of $\mathcal{A}_{\deg i'}$ where $i' = i - \deg(\mathbf{d})$.*

6.5. **An Intermediate Weight Calculation:** $(1 - \mathfrak{t}_{n-1}) r_{\mathrm{BN}}$. In this section we prove that the pattern we noticed in computer-aided calculation for small values of $n$ can be proved to hold for all $n$.

**Theorem 6.8.** *Fix $n \in \mathbb{N}$, and let $K_n = (V, E)$ be the complete graph on vertex set $V = [n]$. Consider the Baker-Norine rank $r_{\mathrm{BN}} \colon \mathbb{Z}^n \to \mathbb{Z}$ on $K_n$. For any $\mathbf{a} \in \mathcal{A}_{\deg i}$,*

$$(37) \quad a_1 + \cdots + a_{n-2} \le i \iff a_n \ge 0 \iff r_{\mathrm{BN}}(\mathbf{a}) = r_{\mathrm{BN}}(\mathbf{a} - \mathbf{e}_{n-1}) + 1.$$

We remark that (37) generalizes Proposition 10 of [CB13].

*Proof.* For all $\mathbf{a} \in \mathcal{A}$, $\mathbf{a} \ge \mathbf{0}$ iff $a_n \ge 0$, since all other coordinates of $\mathbf{a}$ are non-negative. For $\mathbf{a} \in \mathcal{A}_{\deg i}$, in view of (34) when get

$$\mathbf{a} \ge \mathbf{0} \iff a_n \ge 0 \iff a_1 + \cdots + a_{n-2} \le i.$$

Hence Lemma 6.3 implies that for $\mathbf{a} \in \mathcal{A}_{\deg i}$,

$$(38) \quad a_1 + \cdots + a_{n-2} \le i \implies r_{\mathrm{BN}}(\mathbf{a}) = r_{\mathrm{BN}}(\mathbf{a} - \mathbf{e}_{n-1}) + 1.$$

We now prove the reverse implication by, roughly speaking, giving a calculation that shows that there is "no more room" for $r_{\mathrm{BN}}(\mathbf{a}) - r_{\mathrm{BN}}(\mathbf{a} - \mathbf{e}_i)$ to be 1 otherwise, given that we know the offset of $1 + r_{\mathrm{BN}, K_n}$. Let us make this precise.

For any $i \in \mathbb{Z}$, let

$$M_i = \big|\{\mathbf{a} \in \mathcal{A}_{\deg i} \mid r_{\mathrm{BN}}(\mathbf{a}) = r_{\mathrm{BN}}(\mathbf{a} - \mathbf{e}_{n-1}) + 1\}\big|$$

and let

$$N_i = \big|\{\mathbf{a} \in \mathcal{A}_{\deg i} \mid a_1 + \cdots + a_{n-2} \le i\}\big|.$$

Then (38) implies $M_i \ge N_i$, and (37) holds provided that we can show $M_i = N_i$ for all $i$. Since $\mathbf{a} \in \mathcal{A}$ implies that $a_1, \ldots, a_{n-2} \ge 0$, it follows that for $i \le -1$ we have $M_i = N_i = 0$; similarly, since $a_1, \ldots, a_{n-2} \le n - 1$ for $\mathbf{a} \in \mathcal{A}$, we have $a_1 + \cdots + a_{n-2} \le (n-1)(n-2)$; hence for $i \ge n(n-2)$ we have

$$a_1 + \cdots + a_{n-2} \le n(n-2) \le i,$$

and hence for such $i$ we have $N_i = |\mathrm{Pic}_i| = n^{n-2}$, and hence $M_i = n^{n-2}$ as well. Our strategy will be to show that for sufficiently large $\ell \in \mathbb{N}$ we have

$$M_0 + \cdots + M_\ell = N_0 + \cdots + N_\ell;$$

if so, then the inequalities $M_i \geq N_i$ must hold with equality (i.e., there is "no room" for some $N_i$ to be strictly smaller than $M_i$).

Let us take a large $\ell \in \mathbb{N}$; and consider $M_0 + \cdots + M_\ell$: for each $\mathbf{a} \in \mathcal{A}_{\deg \ell}$ we have $r_{\mathrm{BN}}(\mathbf{a}) = \ell - g$ and $r_{\mathrm{BN}}(\mathbf{a} - \mathbf{e}_{n-1}(\ell+1)) = -1$, and hence

(39)
$$\sum_{i=0}^{\ell} (r_{\mathrm{BN}}(\mathbf{a}-i\mathbf{e}_{n-1}) - r_{\mathrm{BN}}(\mathbf{a}-(i+1)\mathbf{e}_{n-1})) = r_{\mathrm{BN}}(\mathbf{a}) - r_{\mathrm{BN}}(\mathbf{a}-\mathbf{e}_{n-1}(\ell+1)) = \ell-g+1.$$

But for all $j$, $\mathcal{A}_j$ is a set of $\mathrm{Pic}_j$ representatives; hence for fixed $i$, as $\mathbf{a}$ varies over $\mathcal{A}_\ell$, and $\mathbf{a} - i\mathbf{e}_n$ varies over a set of $\mathrm{Pic}_{\ell-i}$ representatives; hence

$$\sum_{\mathbf{a} \in \mathcal{A}_\ell} (r_{\mathrm{BN}}(\mathbf{a} - i\mathbf{e}_{n-1}) - r_{\mathrm{BN}}(\mathbf{a} - (i+1)\mathbf{e}_{n-1}))$$
$$= \sum_{p \in \mathrm{Pic}_{\ell-i}} (r_{\mathrm{BN}}(p) - r_{\mathrm{BN}}(p - \mathbf{e}_{n-1}))$$
$$= \sum_{\mathbf{a}' \in \mathcal{A}_{\ell-i}} (r_{\mathrm{BN}}(\mathbf{a}') - r\mathrm{BN}(\mathbf{a}' - \mathbf{e}_{n-1}))$$
$$= M_{\ell-i}$$

(since $r_{\mathrm{BN}}(\mathbf{a}') - r_{\mathrm{BN}}(\mathbf{a}' - \mathbf{e}_{n-1})$ is either 0 or 1, and $M_{\ell-i}$ counts the total number equal to 1). Hence summing (39) over all $\mathbf{a} \in \mathcal{A}_\ell$ we get

(40) $$M_\ell + M_{\ell-1} + \cdots + M_0 = n^{n-2}(\ell - g + 1).$$

Next consider $N_0 + \cdots + N_\ell$ for $\ell$ large: note that for all $(a_1, \ldots, a_{n-2}) \in \{0, \ldots, n-1\}^{n-2}$ and $i \in \mathbb{Z}$, we have

$$\text{either} \quad a_1 + \cdots + a_{n-2} \leq i$$
$$\text{or} \quad a_1 + \cdots + a_{n-2} \geq i + 1$$

(i.e., exactly one of the two inequalities above holds), and hence

$$\text{either} \quad a_1 + \cdots + a_{n-2} \leq i$$
$$\text{or} \quad (n-1-a_1) + \cdots + (n-1-a_{n-2}) \leq (n-1)(n-2) - i - 1.$$

Since $(a_1, \ldots, a_{n-2}) \mapsto (n-1-a_1, \ldots, n-1-a_{n-2})$ is a bijection of $\{0, \ldots, n-1\}^{n-2}$ to itself, it follows that for all $i$ and all $a_1, \ldots, a_{n-2} \in \{0, \ldots, n-1\}$, either $(a_1, \ldots, a_{n-2}) \in \{0, \ldots, n-1\}^{n-2}$ is counted once either in $N_i$, or $(n-1-a_1, \ldots, n-1-a_{n-2})$ is counted once in $N_{(n-2)(n-1)-i-1}$; hence

$$N_i + N_{(n-2)(n-1)-i-1} = n^{n-2}.$$

Hence for all $i \in \mathbb{Z}$ we have

$$N_0 + \cdots + N_{(n-2)(n-1)-1} = \frac{(n-2)(n-1)n^{n-2}}{2},$$

and for $\ell \geq (n-1)(n-2) - 1$ we have

$$N_0 + \ldots + N_\ell = \frac{(n-2)(n-1)n^{n-2}}{2} + n^{n-2}(\ell - (n-1)(n-2) + 1)$$
$$= n^{n-2}\left(\frac{(n-1)(n-2)}{2} + \ell - (n-1)(n-2) + 1\right)$$
$$= n^{n-2}(\ell - g + 1),$$

in view of the fact that
$$g = 1 + |E| - |V| = 1 + \frac{n(n-1)}{2} - n = \frac{2 + n^2 - n - 2n}{2} = \frac{(n-1)(n-2)}{2}.$$

Hence, from (40) we have
$$N_0 + \ldots + N_\ell = n^{n-2}(\ell - g + 1) = M_0 + \cdots + M_\ell$$

for $\ell$ large. But since $M_i \geq N_i$ for all $i$, we must have $N_i = M_i$ for all $0 \leq i \leq \ell$; hence $N_i = M_i$ for all $i$. □

6.6. **A New Rank Formula for the Complete Graph and an Algorithm.** Cori and Le Borgne [CB13] (after Proposition 6, bottom of page 9 and in [CLB16],Proposition 13) describe an $O(n)$ algorithm that computes $r_{\mathrm{BN}}(\mathbf{d})$ for the complete graph $K_n$. Also, they show that when $\mathbf{d}$ is a *sorted parking configuration*, meaning that $0 \leq d_i < i$ for $i < n$ and $d_1 \leq d_2 \leq \cdots \leq d_{n-1}$ (and $d_n$ is unconstrained), they show (see Theorem 12 [CLB16]) that setting
$$q = \lfloor (d_n + 1)/(n-1) \rfloor, \quad r = (d_n + 1) \bmod (n-1)$$

one has
$$r_{\mathrm{BN}}(\mathbf{d}) = -1 + \sum_{i=1}^{n} \max\Big(0, q - i + 1 + d_i + \chi(i \leq r)\Big),$$

where $\chi(P)$ is 1 if $P$ is true, and 0 if $P$ is false.

Here we give another formula for the rank, perhaps related to the above formula; by contrast, our formula holds for $\mathbf{a} \in \mathcal{A}$, but easily generalizes to all $\mathbf{d} \in \mathbb{Z}^n$. The formula is a corollary to Theorem 6.8.

**Corollary 6.9.** *Let $n \in \mathbb{Z}$, and $\mathcal{A}$ be as in (33). For any $\mathbf{a} \in \mathcal{A}$ we have*

$$(41) \quad r_{\mathrm{BN},K_n}(\mathbf{a}) = -1 + \left| \left\{ i = 0, \ldots, \deg(\mathbf{a}) \ \bigg| \ \sum_{j=1}^{n-2} \big((a_j + i) \bmod n\big) \leq \deg(\mathbf{a}) - i \right\} \right|.$$

*In particular, for any $\mathbf{d} \in \mathbb{Z}^n$ we have*
(42)
$$r_{\mathrm{BN},K_n}(\mathbf{d}) = -1 + \left| \left\{ i = 0, \ldots, \deg(\mathbf{d}) \ \bigg| \ \sum_{j=1}^{n-2} \big((d_j - d_{n-1} + i) \bmod n\big) \leq \deg(\mathbf{d}) - i \right\} \right|.$$

*Proof.* Since $\mathbf{a} - (\deg(\mathbf{a}) + 1)\mathbf{e}_{n-1}$ has negative degree, we have

$$(43) \quad \sum_{i=0}^{\deg(\mathbf{a})} \Big( r_{\mathrm{BN}}(\mathbf{a} - i\mathbf{e}_{n-1}) - r_{\mathrm{BN}}(\mathbf{a} - (i+1)\mathbf{e}_{n-1}) \Big) = r_{\mathrm{BN}}(\mathbf{a}) - (-1).$$

According to Theorem 6.8, for a fixed $i$,
$$r_{\mathrm{BN}}(\mathbf{a} - i\mathbf{e}_{n-1}) - r_{\mathrm{BN}}(\mathbf{a} - (i+1)\mathbf{e}_{n-1})$$

equals 1 or 0 according to whether or not the unique $\mathbf{a}' \in \mathcal{A}$ that is equivalent to $\mathbf{a} - i\mathbf{e}_{n-1}$ satisfies

$$(44) \quad a'_1 + \cdots + a'_{n-2} \leq \deg(\mathbf{a}').$$

According to Lemma 6.6, since the $(n-1)$-th component of $\mathbf{a} - i\mathbf{e}_{n-1}$ is $-i$, $\mathbf{a}'$ is given as
$$\forall j \in [n-2], \quad a'_j = (a_j + i) \bmod n,$$

and $(a'_{n-1} = 0)$ and $\deg(\mathbf{a}') = \deg(\mathbf{a}) - i$. Hence (44) holds iff

$$\sum_{j=1}^{n-2}\big((a_j + i) \bmod n\big) \le \deg(\mathbf{a}) - i.$$

Hence, in view of (43) we have (41).

To prove (42), we note that any $\mathbf{d} \in \mathbb{Z}^n$ is equivalent to $\mathbf{a} \in \mathcal{A}$, where

$$a_j = (d_j - d_{n-1}) \bmod n$$

for $j \le n-2$, and $\deg(\mathbf{a}) = \deg(\mathbf{d})$. $\square$

**Remark 6.10.** In the proof above we are making use of the fact that if $f\colon \mathbb{Z}^n \to \mathbb{Z}$ is any function that is initially equal to a constant, then then

$$f(\mathbf{d}) = \Big(\big((1-\mathfrak{t}) + (1-\mathfrak{t}_{n-1})\mathfrak{t}_{n-1} + (1-\mathfrak{t}_{n-1})\mathfrak{t}_{n-1}^2 + \cdots\big)f\Big)(\mathbf{d})$$

where the right-hand-side represents a finite sum, since for any fixed $\mathbf{d}$, for sufficiently large $m \in \mathbb{N}$ we have

$$\big((1-\mathfrak{t}_{n-1})\mathfrak{t}_{n-1}^m f\big)(\mathbf{d}) = 0.$$

One can similarly write, for any $i \in [n]$,

$$(1 - \mathfrak{t}_i)^{-1} = 1 + \mathfrak{t}_i + \mathfrak{t}_i^2 + \cdots$$

with the right-hand-side representing a finite sum when applied to an initially vanishing function $f$ at any given value $\mathbf{d}$. It follows that if $f, f'$ are initially zero, then

(45) $\qquad (1 - \mathfrak{t}_i)f = h \quad \iff \quad f = (1 + \mathfrak{t}_i + \mathfrak{t}_i^2 + \cdots)h.$

At times one of the two conditions above is easier to show that the other, at times not. For example, Theorem 6.8 above gives us a formula for $f = (1 - \mathfrak{t}_{n-1})r_{\mathrm{BN}}$ over $\mathbf{a} \in \mathcal{A}$; in Theorem 6.15 we determine $h = (1 - \mathfrak{t}_n)f$, but it is just as easy to apply either side of (45) with $i = n$. On the other hand, to compute the weight of $r_{\mathrm{BN}}$ in Theorem 6.17, with $h$ as above and

$$W = (1 - \mathfrak{t}_1)\ldots(1 - \mathfrak{t}_{n-2})h,$$

the above formula seems easier to verity than the equivalent

$$h = (1 + \mathfrak{t}_1 + \mathfrak{t}_1^2 + \cdots)\ldots(1 + \mathfrak{t}_{n-2} + \mathfrak{t}_{n-2}^2 + \cdots)W.$$

Next we briefly give a linear time algorithm to compute $r_{\mathrm{BN}}$ of the complete graph based on (41) or (42) in Corollary 6.9.

First, for simplicity, take an arbitrary $\mathbf{d} \in \mathbb{Z}^n$ and note that the equivalent $\mathbf{a} \in \mathcal{A}$ has $a_i = (d_i - d_{n-1}) \bmod n$ for $i \le n-2$ and $\deg(\mathbf{a}) = \deg(\mathbf{d})$. Hence it suffices to show how to compute (41) with $\mathbf{a} \in \mathcal{A}$.

Setting

$$g(i) = \sum_{j=1}^{n-2}\big((a_j + i) \bmod n\big)$$

we have that $g(i+n) = g(i)$ for all $i$, and

(46) $\qquad g(i) = -m_i n + \sum_{j=1}^{n-2} a_j,$

where $m_i$ is the number of $j \in [n-2]$ such that $a_j + i \ge n$, i.e., with $a_j \ge n - i$.

Next, we claim that we can compute $m_0, \ldots, m_{n-1}$ in linear time: indeed, by a single pass through $a_1, \ldots, a_{n-2}$, one can count for each $k = 1, \ldots, n-1$ the number,
$$m'_k = \big|\{j \in [n-2] \mid a_j = k\}\big|,$$
i.e., the number of $j$ for which $a_j = k$; then one computes $m_0, \ldots, m_{n-1}$ by setting $m_0 = 0$ and for $k = 1, \ldots, n-1$ setting $m_k = m'_{n-k} + m_{k-1}$.

Once we compute $m_0, \ldots, m_{n-1}$, we can compute $g(0), \ldots, g(n-1)$ in linear time by computing $\sum_j a_j$ (once) and then applying (46) for each $i = 0, \ldots, n-1$.

Now note that for $k = \{0, \ldots, n-1\}$, we have that for any $i \in \{0, \ldots, \deg(\mathbf{a})\}$ with $i \bmod n = k$, we have $g(i) = g(k)$, and hence the condition
$$\sum_{j=1}^{n-2} \big((a_j + i) \bmod n\big) \le \deg(\mathbf{a}) - i$$
is equivalent to
$$i + g(k) \le \deg(\mathbf{a}),$$
and hence the number of such $i$, for $k$ fixed, is
$$\Big\lfloor \big(\deg(\mathbf{a}) - g(k) + n\big)/n \Big\rfloor.$$
Hence one can write
$$r_{\mathrm{BN}}(\mathbf{a}) = -1 + \sum_{k=0}^{n-1} \Big\lfloor \big(\deg(\mathbf{a}) - g(k) + n\big)/n \Big\rfloor,$$
which completes an $O(n)$ time algorithm to compute $r_{\mathrm{BN}}$.

6.7. **The Second Coordinates for Pic.** To complete our computation of the weight of $r_{\mathrm{BN}}$ of the complete graph, we use a new set of coordinates. As explained in Subsection 6.1, the second coordinates turn out to represent Pic as a product
$$(47) \qquad \mathrm{Pic} = (\mathbb{Z}/n\mathbb{Z})^{n-2} \times \mathbb{Z}.$$

**Notation 6.11.** For any $n \in \mathbb{N}$ and $i \in \mathbb{Z}$, we use
(1) $\mathcal{B} = \mathcal{B}(n)$ to denote the set $\{0, \ldots, n-1\}^{n-2}$ (and usually we just write $\mathcal{B}$ since $n$ will be fixed); and
(2) for any $\mathbf{b} \in \mathcal{B}$ and $i \in \mathbb{Z}$, we use $\langle \mathbf{b}, i \rangle$ to denote
$$(48) \quad \langle \mathbf{b}, i \rangle = (b_1, \ldots, b_{n-2}, 0, i - b_1 - \cdots - b_{n-2}) \in \mathcal{A}_{\deg i} \subset \mathbb{Z}^n_{\deg i} \subset \mathbb{Z}^n.$$
(3) if $\mathbf{c} \in \mathbb{Z}^{n-2}$, we use $\mathbf{c} \bmod n$ to denote the component-wise application of $\bmod n$, i.e.,
$$\mathbf{c} \bmod n = \big(c_1 \bmod n, \ldots, c_{n-2} \bmod n\big) \in \mathcal{B} = \{0, \ldots, n-1\}^{n-2}.$$

**Definition 6.12.** For fixed $n \in \mathbb{Z}$, we refer to $\mathcal{B} = \mathcal{B}(n)$ and the map $\mathcal{B} \times \mathbb{Z} \to \mathbb{Z}^n$ in (48) as the *second coordinates* of $\mathrm{Pic}(K_n)$ representatives.

**Proposition 6.13.** *Let $n \in \mathbb{N}$, and let notation be as in Notation 6.5 and 6.11. Consider the complete graph, $K_n$, and equivalence modulo $\mathrm{Image}(\Delta_{K_n})$. Then:*
*(1) for each $\mathbf{b} \in \mathcal{B}$ and $i \in \mathbb{Z}$,*
$$\langle (b_1, \ldots, b_{n-2}), i \rangle = (a_1, \ldots, a_n),$$
*where*
$$a_1 = b_1, \ \ldots, \ a_{n-2} = b_{n-2}, \ a_{n-1} = 0,$$

*and*
$$a_n = i - b_1 - \cdots - b_{n-2}.$$

(2) *For all $i \in \mathbb{Z}$, the set $\mathcal{B} \times \{i\}$ is taken via $\langle \cdot, \cdot \rangle$ bijectively to $\mathcal{A}_{\deg i}$, and hence to a set of representatives of $\mathrm{Pic}_i$.*

(3) *For all $i \in \mathbb{Z}$, each $\mathbf{d} \in \mathbb{Z}_{\deg i}^n$ is equivalent to a unique element of the form $\langle \mathbf{b}, i \rangle$ with $\mathbf{b} \in \mathcal{B}$, namely with*
$$\mathbf{b} = (d_1 - d_{n-1}, \ldots, d_{n-2} - d_{n-1}) \bmod n,$$
*where $\bmod\ n$ is the component-wise application of $\bmod\ n$, i.e., $b_i = (d_i - d_{n-1}) \bmod n \in \{0, \ldots, n-1\}$.*

(4) *For any $\mathbf{b}, \mathbf{b}' \in \mathcal{B} = \{0, \ldots, n-1\}^{n-2}$ and any $i, i' \in \mathbb{Z}$, we have*
$$\langle \mathbf{b}, i \rangle + \langle \mathbf{b}', i' \rangle \sim \langle (\mathbf{b} + \mathbf{b}') \bmod n, i + i' \rangle.$$
*Similarly for subtraction, i.e., with $-$ everywhere replacing $+$.*

*Proof.* (1) is immediate from the notation. (2) follows from (1). (3) follows from (1) and Lemma 6.6. (4) follows from (3). □

**Example 6.14.** Applying the above proposition, we see that

(49)     $\mathbf{e}_1 \sim \langle \mathbf{e}_1, 1 \rangle, \ \ldots, \mathbf{e}_{n-2} \sim \langle \mathbf{e}_{n-2}, 1 \rangle, \ \mathbf{e}_{n-1} \sim \langle (n-1)\mathbf{1}, 1 \rangle, \ \mathbf{e}_n \sim \langle \mathbf{0}, 1 \rangle,$

where we use $\mathbf{e}_i$ to denote the vector in $\mathbb{Z}^n$ or in $\mathbb{Z}^{n-2}$, as appropriate. Moreover, equality holds in all the above, except for $\mathbf{e}_{n-1}$, where
$$\mathbf{e}_{n-1} \sim \langle (n-1)\mathbf{1}, 1 \rangle = (n-1, \ldots, n-1, 0, 1-(n-2)(n-1)).$$

6.8. **Computation of $(1 - \mathfrak{t}_n)(1 - \mathfrak{t}_{n-1})r_{\mathrm{BN}}$.**

**Theorem 6.15.** *Fix $n \in \mathbb{N}$, and let $K_n = (V, G)$ be the complete graph on vertex set $V = [n]$, i.e., $E$ consists of exactly one edge joining any two distinct vertices. Consider the Baker-Norine rank $r_{\mathrm{BN}} \colon \mathbb{Z}^n \to \mathbb{Z}$ on $K_n$.*

*(1) If $\mathbf{a} \in \mathcal{A}_{\deg i}$, then*

(50)     $(1 - \mathfrak{t}_n)(1 - \mathfrak{t}_{n-1})r_{\mathrm{BN},K_n}(\mathbf{a}) = \begin{cases} 1 & \text{if } a_1 + \cdots + a_{n-2} = i, \text{ and} \\ 0 & \text{otherwise.} \end{cases}$

*(2) For all $\mathbf{b} \in \mathcal{B}$ and $i \in \mathbb{Z}$,*

(51)     $(1 - \mathfrak{t}_n)(1 - \mathfrak{t}_{n-1})r_{\mathrm{BN},K_n}(\langle \mathbf{b}, i \rangle) = \begin{cases} 1 & \text{if } b_1 + \cdots + b_{n-2} = i, \text{ and} \\ 0 & \text{otherwise.} \end{cases}$

*Proof.* The left-hand-side of (50) equals
$$(1 - \mathfrak{t}_n)(1 - \mathfrak{t}_{n-1})r_{\mathrm{BN},K_n}(\mathbf{a}) = (1 - \mathfrak{t}_{n-1})r_{\mathrm{BN},K_n}(\mathbf{a}) - (1 - \mathfrak{t}_{n-1})r_{\mathrm{BN},K_n}(\mathbf{a} - \mathbf{e}_n).$$
Note that if $\mathbf{a} \in \mathcal{A}_{\deg i}$, then
$$\mathbf{a} - \mathbf{e}_n = (a_1, \ldots, a_{n-2}, 0, i - 1 - a_1 - \cdots - a_{n-2}) \in \mathcal{A}_{\deg i - 1}.$$
By Theorem 6.8, $(1 - \mathfrak{t}_{n-1})r_{\mathrm{BN},K_n}(\mathbf{a})$ is 1 or 0 according to whether or not $a_1 + \cdots + a_{n-2} \le i$ or not, and similarly with $\mathbf{a}$ replaced by $\mathbf{a} - \mathbf{e}_n \in \mathcal{A}_{\deg i - 1}$, according to whether or not $a_1 + \cdots + a_{n-2} \le i - 1$. Hence we conclude (50).

(2) (i.e., (51)) follows immediately from (1) (i.e., (50)). □

When going through the weight calculations in the next two sections, it may be helpful to visualize consequences of Theorem 6.8 in the case $n = 4$, and to consider what (51) means in terms of the $\langle \mathbf{b}, i \rangle$ coordinates, namely that $b_1 + b_2 = i$; see Figure 1.

| $i=0$ | | | | | $i=1$ | | | | | $i=2$ | | | | |
|---|---|---|---|---|---|---|---|---|---|---|---|---|---|---|
|  | 0 | 1 | 2 | 3 |  | 0 | 1 | 2 | 3 |  | 0 | 1 | 2 | 3 |
| 0 | 1 | 0 | 0 | 0 | 0 | 0 | 1 | 0 | 0 | 0 | 0 | 0 | 1 | 0 |
| 1 | 0 | 0 | 0 | 0 | 1 | 1 | 0 | 0 | 0 | 1 | 0 | 1 | 0 | 0 |
| 2 | 0 | 0 | 0 | 0 | 2 | 0 | 0 | 0 | 0 | 2 | 1 | 0 | 0 | 0 |
| 3 | 0 | 0 | 0 | 0 | 3 | 0 | 0 | 0 | 0 | 3 | 0 | 0 | 0 | 0 |

| $i=3$ | | | | | $i=4$ | | | | | $i=5$ | | | | |
|---|---|---|---|---|---|---|---|---|---|---|---|---|---|---|
|  | 0 | 1 | 2 | 3 |  | 0 | 1 | 2 | 3 |  | 0 | 1 | 2 | 3 |
| 0 | 0 | 0 | 0 | 1 | 0 | 0 | 0 | 0 | 0 | 0 | 0 | 0 | 0 | 0 |
| 1 | 0 | 0 | 1 | 0 | 1 | 0 | 0 | 0 | 1 | 1 | 0 | 0 | 0 | 0 |
| 2 | 0 | 1 | 0 | 0 | 2 | 0 | 0 | 1 | 0 | 2 | 0 | 0 | 0 | 1 |
| 3 | 1 | 0 | 0 | 0 | 3 | 0 | 1 | 0 | 0 | 3 | 0 | 0 | 1 | 0 |

| $i=6$ | | | | |
|---|---|---|---|---|
|  | 0 | 1 | 2 | 3 |
| 0 | 0 | 0 | 0 | 0 |
| 1 | 0 | 0 | 0 | 0 |
| 2 | 0 | 0 | 0 | 0 |
| 3 | 0 | 0 | 0 | 1 |

FIGURE 1. The non-zero values of of $(1-\mathfrak{t}_{n-1})(1-\mathfrak{t}_n)r_{\mathrm{BN}}(\langle b, i\rangle)$ for $n=4$, $\mathbf{b}=(b_1,b_2)\in\{0,1,2,3\}^2$, namely 1 if $b_1+b_2=i$, and 0 otherwise.

### 6.9. A Generalization of the Weight Calculation.
To compute the weight of the Baker-Norine rank on $K_n$, we need to apply
$$(1-\mathfrak{t}_1)\ldots(1-\mathfrak{t}_{n-2}).$$
However, (51) implies that
$$(1-\mathfrak{t}_n)(1-\mathfrak{t}_{n-1})r_{\mathrm{BN},K_n}(\langle\mathbf{b},i\rangle)=g(b_1+\cdots+b_{n-2}-i),$$
for some function $g$ (namely the "Dirac delta function at 0," i.e., the function that is 1 at 0 and otherwise 0). We find it conceptually simpler to prove a theorem that applies
$$(1-\mathfrak{t}_1)\ldots(1-\mathfrak{t}_{n-2})$$
to any function of $\langle b, i\rangle$ of the form
$$g(b_1+\cdots+b_{n-2}-i).$$
Here is the result.

It will be helpful to introduce the following "tensor" notation: if $J\subset[n-2]$, then set
$$\mathfrak{t}_J=\prod_{j\in J}\mathfrak{t}_j. \tag{52}$$

**Proposition 6.16.** *Let $h\colon\mathbb{Z}^n\to\mathbb{Z}$ be any function that is invariant under translation by the image of the Laplacian of the complete graph. Say that for all $(\mathbf{b},i)\in\mathcal{B}\times\mathbb{Z}$, $h(\langle\mathbf{b},i\rangle)=g(b_1+\cdots+b_{n-2}-i)$ for some function $g$, i.e., $h$ depends only on the value of $b_1+\cdots+b_{n-2}-i$. Then*

(1) if $j \in [n-2]$ and $\mathbf{b} \in \mathcal{B} = \{0, \cdots, n-1\}^{n-2}$ has $b_j > 0$, then for all $i \in \mathbb{Z}$ we have

(53) $$((1 - \mathfrak{t}_j)h)(\langle b, i \rangle) = 0;$$

(2) let $j \in [n-2]$ and $J' \subset [n-2]$ with $j \notin J'$; if $\mathbf{b} \in \mathcal{B} = \{0, \cdots, n-1\}^{n-2}$ has $b_j > 0$, then for all $i \in \mathbb{Z}$ we have

(54) $$\big((1 - \mathfrak{t}_j)\mathfrak{t}_{J'}h\big)(\langle b, i \rangle) = 0$$

(using the "tensor" notation (52));

(3) if $\mathbf{b} \in \mathcal{B}$ with $\mathbf{b} \ne \mathbf{0}$ (hence $b_j > 0$ for some $j \in [n-2]$),

(55) $$\big((1 - \mathfrak{t}_1) \ldots (1 - \mathfrak{t}_{n-2})h\big)(\langle \mathbf{b}, i \rangle) = 0;$$

and

(4) (in the remaining case, $\mathbf{b} = \mathbf{0}$)

(56) $$\big((1 - \mathfrak{t}_1) \ldots (1 - \mathfrak{t}_{n-2})h\big)(\langle \mathbf{0}, i \rangle) = \sum_{k=0}^{n-2} (-1)^k \binom{n-2}{k} g(i - kn).$$

We remark that the proof below shows that claims (1) and (2) above hold, more generally, whenever

$$h(\langle \mathbf{b}, i \rangle) = g(b_1, \ldots, b_{j-1}, b_j - i, b_{j+1}, \ldots, b_{n-2})$$

for some $g$, i.e., $h$ is an arbitrary function, except that its dependence on $b_j$ and $i$ is only on $b_j - i$ and the rest of the $b_{j'}$ with $j' \ne j$.

*Proof.* Our proof will constantly use (49).

Proof of (1): if $b_j > 0$, then $\mathbf{b} - \mathbf{e}_j \in \mathcal{B}$, and hence

$$\langle \mathbf{b}, i \rangle - \mathbf{e}_j = \langle \mathbf{b} - \mathbf{e}_j, i - 1 \rangle,$$

and hence

$$\big((1 - \mathfrak{t}_j)h\big)(\langle b, i \rangle) = h(\langle \mathbf{b}, i \rangle) - h(\langle \mathbf{b} - \mathbf{e}_j, i - 1 \rangle)$$
$$= g\big((b_1 + \cdots + b_{n-2}) - i\big) - g\big((b_1 + \cdots + b_{n-2} - 1) - (i - 1)\big) = 0.$$

This gives (53).

Proof of (2): let

$$\mathbf{b}' = (\mathbf{b} - \mathbf{e}_{J'}) \bmod n.$$

Since $j \notin J'$ we have $b'_j = b_j > 0$, and hence $\mathbf{b}' - \mathbf{e}_j \in \mathcal{B}$. Hence

$$\big(\mathfrak{t}_J h\big)(\langle \mathbf{b}, i \rangle) = h\big(\langle \mathbf{b}', i - |J'| \rangle\big)$$
$$\big(\mathfrak{t}_j \mathfrak{t}_{J'} h\big)(\langle \mathbf{b}, i \rangle) = h\big(\langle \mathbf{b}' - \mathbf{e}_j, i - |J'| - 1 \rangle\big).$$

Hence the same calculation as in the previous paragraph (with $\mathbf{b}'$ replacing $\mathbf{b}$ and $i - |J'|$ replacing $i$) gives (54).

Proof of (3): we have

$$(1 - \mathfrak{t}_1) \ldots (1 - \mathfrak{t}_{n-2}) = \sum_{J' \subset [n-2] \setminus \{j\}} (-1)^{|J'|} (1 - \mathfrak{t}_j) \mathfrak{t}_{J'},$$

and so (54) implies (55).

Proof of (4): for any $J \subset [n-2]$, using (49) we have

$$\langle \mathbf{0}, i \rangle - \mathbf{e}_J \sim \langle (n-1)\mathbf{e}_J, i - |J| \rangle,$$

and hence
$$f\bigl(\langle \mathbf{0}, i\rangle - \mathbf{e}_J\bigr) = f\bigl(\langle (n-1)\mathbf{e}_J, i - |J|\rangle\bigr) = g((n-1)|J| - i + |J|) = g(n|J| - i).$$
Since
$$(1 - \mathfrak{t}_1)\ldots(1 - \mathfrak{t}_{n-2}) = \sum_{J\subset [n-2]} (-1)^{|J|}\mathfrak{t}_J,$$
we get
$$\bigl((1 - \mathfrak{t}_1)\ldots(1 - \mathfrak{t}_{n-2})h\bigr)(\langle \mathbf{0}, i\rangle) = \sum_{J\subset [n-2]} (-1)^{|J|} g\bigl(n|J| - i\bigr)$$
and (56) follows. □

### 6.10. Computation of $W$.

**Theorem 6.17.** *Fix $n \in \mathbb{N}$, and let $K_n = (V, E)$ be the complete graph on vertex set $V = [n]$. Consider the Baker-Norine rank $r_{\mathrm{BN}}\colon \mathbb{Z}^n \to \mathbb{Z}$ on $K_n$. The weight, $W = \mathfrak{m}(r_{\mathrm{BN}, K_n})$, is given by*
(57)
$$W(\langle \mathbf{b}, i\rangle) = \begin{cases} (-1)^\ell \binom{n-2}{\ell} & \text{if } \mathbf{b} = \mathbf{0} \text{ and } i = n\ell \text{ for some } \ell = 0, \ldots, n-2, \text{ and} \\ 0 & \text{otherwise.} \end{cases}$$

*Proof.* Setting
$$h(\langle \mathbf{b}, i\rangle) = \bigl((1 - \mathfrak{t}_{n-1})(1 - \mathfrak{t}_n)r_{\mathrm{BN}}\bigr)(\langle \mathbf{b}, i\rangle),$$
(51) shows that
$$h(\langle \mathbf{b}, i\rangle) = g(b_1 + \cdots + b_{n-2} - i),$$
where $g(0) = 1$ and elsewhere $g$ vanishes. Since
$$W = (1 - \mathfrak{t}_1)\cdots(1 - \mathfrak{t}_{n-2})h,$$
we may apply Proposition 6.16 and conclude: (1) if $\mathbf{b} \in \mathcal{B}$ is nonzero, then (55) implies that
$$W(\langle \mathbf{b}, i\rangle) = 0,$$
and (2) if $\mathbf{b} = \mathbf{0}$, then
$$W(\langle \mathbf{0}, i\rangle) = \sum_{k=0}^{n-2} (-1)^k \binom{n-2}{k} g(nk - i).$$
Hence $W(\langle \mathbf{0}, i\rangle) = 0$ unless $i$ is of the form $nk$, with $0 \le k \le n-2$, in which case
$$W(\langle \mathbf{0}, nk\rangle) = (-1)^k \binom{n-2}{k}.$$
□

### 6.11. Remark on Theorem 6.17.
Another important consequence of Theorem 6.17 is that, by symmetry, for any $\mathbf{d} \in \mathbb{Z}^n$, and any distinct $i, j \in [n]$ we have
$$\bigl((1 - \mathfrak{t}_i)(1 - \mathfrak{t}_j)W\bigr)(\mathbf{d}) \ge 0.$$
In [FF] this will imply that when we can model $f = 1 + r_{\mathrm{BN}, K_n}$ as Euler characteristics of a family of sheaves in a sense explained there.

## 7. Fundamental Domains and the Proofs of Theorems 3.3 and 3.4

In this section we prove the Theorems 3.3 and 3.4. We do so with a tool that we call a *cubism* of $\mathbb{Z}^n$. However, Theorems 3.3 has a more direct proof without using cubisms, so we first give the direct proof. In fact, the direct proof will motivate the definition of a cubism.

### 7.1. Proof of Theorem 3.3 Without Reference to Cubisms.

**Lemma 7.1.** *Let $n \in \mathbb{Z}$, and let $\mathcal{D}^n_{\mathrm{coord}} \subset \mathbb{Z}^n$ given by*

(58) $$\mathcal{D}^n_{\mathrm{coord}} = \{\mathbf{d} \ |d_i = 0 \text{ for at least one } i \in [n]\}.$$

*Then for any $f \colon \mathcal{D}^n_{\mathrm{coord}} \to \mathbb{Z}$, there exist functions $h_i \colon \mathbb{Z}^n \to \mathbb{Z}$ for each $i \in [n]$ such that*

(1) $h_i = h_i(\mathbf{d})$ *is independent of the $i$-th variable, $d_i$, and*
(2)

(59) $$\forall \mathbf{d} \in \mathcal{D}^n_{\mathrm{coord}}, \quad f(\mathbf{d}) = \sum_{i=1}^n h_i(\mathbf{d}).$$

Hence the function $\sum_i h_i$ above is an extension of $f$ to all of $\mathbb{Z}^n$ such that each $h_i$ is independent of its $i$-th variable.

Before giving the formal proof, let us explain the ideas for small $n$. The case $n = 1$ is immediate. The proof for $n = 2$ is as follows: consider

(60) $$g(d_1, d_2) = f(d_1, 0) + f(0, d_2) - f(0, 0):$$

since

$$g(d_1, 0) = f(d_1, 0) + f(0, 0) - f(0, 0) = f(d_1, 0)$$

we have $f(\mathbf{d}) = g(\mathbf{d})$ whenever $d_2 = 0$; by symmetry, the same is true if $d_1 = 0$; hence $g = f$ on all of $\mathcal{D}^2_{\mathrm{coord}}$. But we easily write the right-hand-side of (60) as $h_1(d_2) + h_2(d_1)$, by setting, say, $h_1(d_2) = f(0, d_2) - f(0, 0)$ and setting $h_2(d_1) = f(d_1, 0)$.

Similarly for $n = 3$, and

$$g(d_1, d_2, d_3) = f(d_1, d_2, 0) + f(d_1, 0, d_3) + f(0, d_2, d_3) - f(d_1, 0, 0) - f(0, d_2, 0) - f(0, 0, d_3) + f(0, 0, 0).$$

For all $n \geq 4$, we simply need to introduce convenient notation.

*Proof.* For $\mathbf{d} \in \mathbb{Z}^n$ and $I \subset [n]$, introduce the notation

$$\mathbf{d}_I = \sum_{i \in I} d_i \mathbf{e}_i.$$

Consider the function $g \colon \mathbb{Z}^n \to \mathbb{Z}$ given by

(61) $$g(\mathbf{d}) = \sum_{I \subset [n], \ I \neq [n]} f(\mathbf{d}_I)(-1)^{n-1-|I|}$$

(which makes sense, since $\mathbf{d}_I \in \mathcal{D}^n_{\mathrm{coord}}$ whenever $I \neq [n]$). We claim that $g = f$ when restricted to $\mathbf{d} \in \mathcal{D}^n_{\mathrm{coord}}$; by symmetry it suffices to check the case $d_n = 0$, whereupon the term $f(\mathbf{d}_I)$ with $n \notin I$ cancels the term corresponding to $I \cup n$, except for the single remaining term where $I = \{1, \ldots, n-1\}$. Hence for $d_n = 0$, $g(\mathbf{d}) = f(\mathbf{d})$, and, by symmetry, $g = f$ on all of $\mathcal{D}^n_{\mathrm{coord}}$.

Now we see that the right-hand-side (61) is of the desired form $\sum_i h_i$ as in the statement of the lemma, by setting

$$h_i = \sum_{i \notin I,\ 1,\ldots,i-1 \in I} f(\mathbf{d}_I)(-1)^{n-1-|I|};$$

since for each $I \subset [n]$ with $I \neq [n]$ there is a unique $i \in [n]$ such that $i \notin I$ but $1, \ldots, i-1 \in I$ (namely the lowest value of $i$ not in $I$), we have $\sum_i h_i$ equals the right-hand-side (61). $\square$

**Theorem 7.2.** *Let $n \in \mathbb{N}$ and $\mathcal{D}^n_{\mathrm{coord}}$ be as in (58). Then any function $f \colon \mathcal{D}^n_{\mathrm{coord}} \to \mathbb{Z}$ has a unique extension to a modular function $h \colon \mathbb{Z}^n \to \mathbb{Z}$.*

*Proof.* The existence of the extension of $h$ is guaranteed by Lemma 7.1. Let us prove uniqueness. By symmetry it suffices to show that the values of $h$ on the set

$$\mathbb{N}^n = \{\mathbf{d} \mid d_i > 0 \text{ for all } i \in [n]\}$$

are uniquely determined. But if $h$ is modular, then

(62) $$h(\mathbf{d}) = \sum_{I \subset [n],\ I \neq \emptyset} (-1)^{|I|+1} h(\mathbf{d} - \mathbf{e}_I).$$

Now we prove by induction on $m$ that for all $m \geq n$, if $\mathbf{d} \in \mathbb{N}^n$ and $\deg(\mathbf{d}) = m$, then $h(\mathbf{d})$ is uniquely determined. The base case is $m = n$, where the only element of degree $n$ in $\mathbb{N}^n$ is $\mathbf{d} = \mathbf{1}$. But for each $I \subset [n]$ with $I \neq \emptyset$, $\mathbf{1} - \mathbf{e}_I \in \mathcal{D}^n_{\mathrm{coord}}$; hence (62) uniquely determines $h(\mathbf{1})$. To prove the inductive claim: let $\mathbf{d} \in \mathbb{N}^n$ with $\deg(\mathbf{d}) = m$; for all $I \subset [n]$ with $I \neq \emptyset$, $\mathbf{d} - \mathbf{e}_I \geq \mathbf{0}$ and $\mathbf{d} - \mathbf{e}_I$ and has degree less than $m$. Hence (62) determines $h(\mathbf{d})$ in terms of values of $h$ that, by induction, have already been determined. $\square$

*Proof of Theorem 3.3.* One direction is immediate; it suffices to show that any modular function, $h$, can be written as a sum of functions, each of which depends on only $n-1$ of its variables. So consider the restriction of $h$ to $\mathcal{D}^n_{\mathrm{coord}}$; then this restriction determines a unique modular function, which must be $h$. But then Theorem 7.2 implies that $h = \sum_i h_i$, where each $h_i$ is independent of its $i$-th variable. $\square$

### 7.2. Fundamental Modular Domains.
Let us restate what we proved in the previous subsection.

**Definition 7.3.** Let $\mathcal{D} \subset \mathbb{Z}^n$. We call $\mathcal{D}$ a *fundamental modular domain* (respectively *subfundamental, superfundamental*) if for every function $f \colon \mathcal{D} \to \mathbb{Z}$ there exists a unique (respectively, at least one, at most one) modular function $h \colon \mathbb{Z}^n \to \mathbb{Z}$ such that $f = h$ on $\mathcal{D}$.

We remark that our terminology results from the following almost immediate facts: a subset of a subfundamental modular domain is subfundamental, and a strict subset of a fundamental domain is not fundamental; similarly for supersets and superfundamental domains.

In the last subsection, Theorem 3.3 was proven via Theorem 7.2, which proved that $\mathcal{D}^n_{\mathrm{coord}}$ is a fundamental modular domain. Theorem 3.4 essentially states that for any $n \in \mathbb{N}$ and $a \in \mathbb{Z}$,

$$\mathcal{D} = \{\mathbf{d} \in \mathbb{Z}^n \mid a \leq \deg(\mathbf{d}) \leq a + n - 1\}$$

is a fundamental modular domain. We can prove both ideas by the method of a *cubism*, that we now explain.

**7.3. Cubisms: Motivation, Definition, and Implication of Domain Fundamentality.** The proof of Theorem 7.2 can be viewed as follows: we ordered the elements of $\mathbb{N}^n$ by a function

$$\operatorname{rank}(\mathbf{d}) = d_1 + \cdots + d_n - (n-1),$$

(so the minimum rank of an element of $\mathbb{Z}^n$ is 1), and proved by induction on $m \geq 1$ that there is a unique extension of a function $h \colon \mathcal{D}^n_{\text{coord}} \to \mathbb{Z}$ to all points of rank at most $m$ so that $(\mathfrak{m}h)(\mathbf{d}) = 0$ for all $\mathbf{d}$ of rank at most $m$. Let us generalize this idea.

**Definition 7.4.** For $\mathbf{d} \in \mathbb{Z}^n$, the $\mathbf{d}$-*cube* refers to the set

$$\operatorname{Cube}(\mathbf{d}) = \{\mathbf{d}' \in \mathbb{Z}^n \mid \mathbf{d} - \mathbf{1} \leq \mathbf{d}' \leq \mathbf{d}\}.$$

We refer to the set of all $\mathbf{d}$-cubes as the set of $n$-*cubes*. If $\mathcal{D} \subset \mathbb{Z}^n$, we say that function $r \colon \mathbb{Z}^n \to \mathbb{N}$ is a *cubism of* $\mathcal{D}$ if, setting

$$\mathcal{D}_m = \mathcal{D} \cup \bigcup_{r(\mathbf{d}) \leq m} \operatorname{Cube}(\mathbf{d}) \tag{63}$$

for $m \in \mathbb{Z}_{\geq 0}$ (hence $\mathcal{D}_0 = \mathcal{D}$), we have

(1) if $m \geq 1$ and $r(\mathbf{d}) = r(\mathbf{d}') = m$, then

$$\operatorname{Cube}(\mathbf{d}) \cap \operatorname{Cube}(\mathbf{d}') \in \mathcal{D}_{m-1}, \tag{64}$$

and
(2) for all $m \geq 1$ and $\mathbf{d} \in \mathbb{Z}^n$ with $r(\mathbf{d}) = m$ we have

$$\bigl|\operatorname{Cube}(\mathbf{d}) \setminus \mathcal{D}_{m-1}\bigr| = 1. \tag{65}$$

In the last paragraph of this section we remark that in some cubisms it is more convenient to replace the partial ordering of the $n$-cubes induced by the function $r \colon \mathbb{Z}^n \to \mathbb{N}$ above with, more generally, a well-ordering or a partial ordering such that each subset has a minimal element.

**Example 7.5.** In Figure 2 we illustrate an example of a cubism of $\mathcal{D}$, with $\mathcal{D} = \mathcal{D}_{\text{coord}}$ as above, suggested by the above proof of Theorem 7.2 and $n = 2$ (so the $n$-cubes are really squares).

**Proposition 7.6.** *If* $\mathcal{D} \subset \mathbb{Z}^n$ *has a cubism, then* $\mathcal{D}$ *is fundamental.*

*Proof.* Fix a function $f \colon \mathcal{D} \to \mathbb{Z}$, and set $g_0 = f$.

Let us prove by induction on $m \in \mathbb{N}$ that there is a unique function $\mathcal{D}_m \to \mathbb{Z}$ such that

(1) $(\mathfrak{m}g_m)(\mathbf{d}) = 0$ for all $\mathbf{d}$ with $r(\mathbf{d}) \leq m$;
(2) the restriction of $g_m$ to $\mathcal{D}_{m-1}$ equals $g_{m-1}$; and
(3) the value of $g_m$ on each $\mathbf{c} \in \mathcal{D}_m \setminus \mathcal{D}_{m-1}$ is determined by the equation $(\mathfrak{m}g_m)(\mathbf{d}) = 0$ for a unique $\mathbf{d} \in \mathcal{D}_{m-1}$ such that $\mathbf{c} \in \operatorname{Cube}(\mathbf{d}) \setminus \mathcal{D}_{m-1}$, via the equation

$$-g_m(\mathbf{c})(-1)^{\deg(\mathbf{d}-\mathbf{c})} = \sum_{\mathbf{c}' \in \operatorname{Cube}(\mathbf{d}) \setminus \{\mathbf{c}\}} g_{m-1}(\mathbf{c}')(-1)^{\deg(\mathbf{d}-\mathbf{c}')}. \tag{66}$$

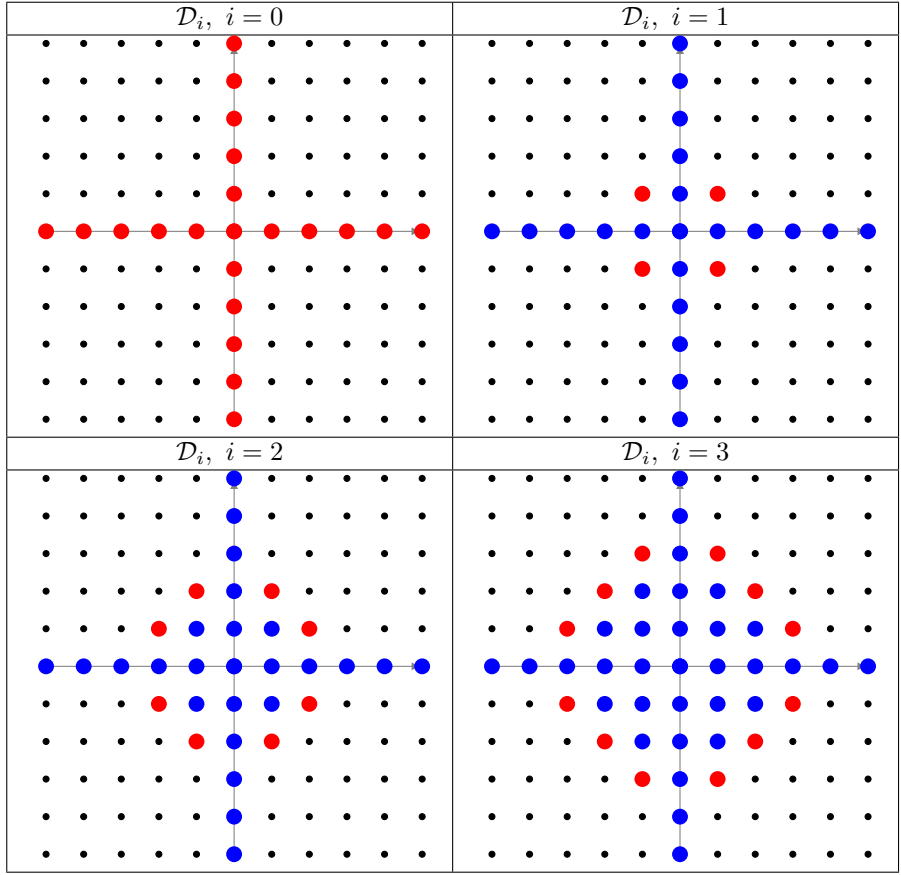

(A) New points $\mathcal{D}_i \setminus \mathcal{D}_{i-1}$ in red, old points, $\mathcal{D}_{i-i}$ in blue

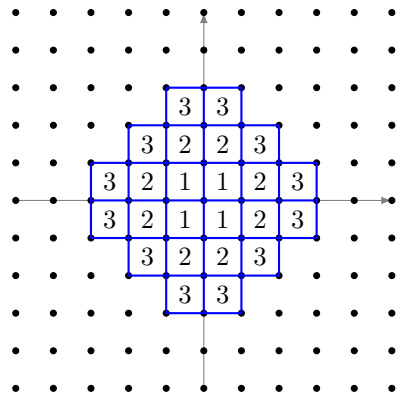

(B) The cubism after 4 steps.

FIGURE 2. A cubism for $\mathcal{D}_{\mathrm{coord}}^n$ with $n = 2$.

The base case $m = 1$ is argued almost exactly as the inductive claim from $m - 1$ to $m$; so we will prove the base case $m = 1$, leaving in $m$ everywhere.

For $m = 1$, we have that $\mathcal{D}_{m-1} = \mathcal{D}_0 = \mathcal{D}$, and (65) implies that for each $\mathbf{d}$ with $r(\mathbf{d}) = m$, there is a unique $\tilde{\mathbf{d}} \notin \mathcal{D}_{m-1}$ in $\mathrm{Cube}(\mathbf{d})$; the equation $(\mathfrak{m}g)(\mathbf{d}) = 0$ is equivalent to

$$(67) \qquad \sum_{\mathbf{c} \in \mathrm{Cube}(\mathbf{d})} g_m(\mathbf{c})(-1)^{\deg(\mathbf{d}-\mathbf{c})} = 0.$$

This determines $g_m(\tilde{\mathbf{d}})$ via (66) with $\mathbf{c} = \tilde{\mathbf{d}}$, since all other $\mathbf{c} \in \mathrm{Cube}(\mathbf{d})$ in the sum (67) either lie in $\mathcal{D}$ or have rank at most $m - 1$; (64) shows that for distinct $\mathbf{d}, \mathbf{d}'$ of rank $m$, the corresponding $\tilde{\mathbf{d}}, \tilde{\mathbf{d}}'$ are distinct, so that it is possible to set the value of $g_m$ as required on all $\tilde{\mathbf{d}}$ that are the unique element of $\mathrm{Cube}(\mathbf{d}) \setminus \mathcal{D}_{m-1}$ for some $\mathbf{d}$ of rank $m$.

For the inductive step, we assume the claim holds for $m - 1$, and we repeat the same argument above. This shows that $g_m \colon \mathcal{D}_m \to \mathbb{Z}$ exist for all $m$ with the desired properties.

Now define $h \colon \mathbb{Z}^n \to \mathbb{Z}$ as follows: for any $\mathbf{d} \in \mathbb{Z}^n$, we have $\mathbf{d} \in \mathrm{Cube}(\mathbf{d}) \subset \mathcal{D}_m$, where $m = r(\mathbf{d})$; hence $g_m(\mathbf{d})$ is defined; set $h(\mathbf{d}) = g_m(\mathbf{d})$.

We claim that $h$ above is modular: indeed, for any $\mathbf{d} \in \mathbb{Z}^n$, if $m = r(\mathbf{d})$, then $\mathfrak{m}g_m(\mathbf{d}) = 0$ and $\mathcal{D}_m$ contains $\mathrm{Cube}(\mathbf{d})$; since $g_{m+1}, g_{m+2}, \ldots$ are all extensions of $g_m$, we have $\mathfrak{m}h(\mathbf{d}) = \mathfrak{m}g_m(\mathbf{d}) = 0$.

Now we claim that $h$ is the unique modular function $\mathbb{Z}^n \to \mathbb{Z}$ whose restriction to $\mathcal{D}$ is $f$: indeed, assume that $h'$ is another such modular function, and that $h \neq h'$; then the definition of $h$ implies that there exists an $m$ such that $g_m$ does not equal the restriction of $h'$ to $\mathcal{D}_m$; consider the smallest such $m$. Since the restrictions of $h$ and $h'$ to $\mathcal{D}_0 = \mathcal{D}$ both equal $f$, we must have $m \geq 1$. It follows that $h(\mathbf{c}) \neq h'(\mathbf{c})$ for some $\mathbf{c} \in \mathcal{D}_m \setminus \mathcal{D}_{m-1}$ with $m \geq 1$; fix such a $\mathbf{c}$. By condition (3) on $g_m$ (i.e., (66) and above), there is some $\mathbf{d}$ with $r(\mathbf{d}) = m$ for which $\mathbf{c}$ which is the unique element of $\mathrm{Cube}(\mathbf{d}) \setminus \mathcal{D}_{m-1}$. But since $h, h'$ agree on $g_{m-1}$, we have

$$(\mathfrak{m}h')(\mathbf{d}) = h'(\mathbf{c})(-1)^{\deg(\mathbf{d}-\mathbf{c})} + \sum_{\mathbf{c}' \in \mathrm{Cube}(\mathbf{d}) \setminus \{\mathbf{c}\}} g_{m-1}(\mathbf{c}')(-1)^{\deg(\mathbf{d}-\mathbf{c}')}$$

$$\neq h(\mathbf{c})(-1)^{\deg(\mathbf{d}-\mathbf{c})} + \sum_{\mathbf{c}' \in \mathrm{Cube}(\mathbf{d}) \setminus \{\mathbf{c}\}} g_{m-1}(\mathbf{c}')(-1)^{\deg(\mathbf{d}-\mathbf{c}')} = 0,$$

and hence $(\mathfrak{m}h')(\mathbf{d}) \neq 0$; hence $h'$ is not modular. $\square$

[Straying a bit, one could define a subcubism by replacing the $= 1$ in (65) by $\geq 1$, and the same proof shows that a $\mathcal{D}$ with a subcubism is subfundamental; similarly for supercubism and $\leq 1$.]

7.4. **Second Proof of Theorem 7.2.** The proof of Theorem 7.2 above can be viewed as giving a cubism (e.g., Figure 2 for $n = 2$). Let us formalize this.

*Second proof of Theorem 7.2.* For each $\mathbf{d} \in \mathbb{Z}^n$, let

$$r(\mathbf{d}) = |d_1| + \cdots + |d_n| + \big|\{i \in [n] \mid d_i \leq 0\}\big| - n + 1;$$

more intuitively, $r(\mathbf{d})$ is just the $L^1$ distance of the furthest point in $\mathrm{Cube}(\mathbf{d})$ to $\mathcal{D}^n_{\mathrm{coord}}$, since if all $d_i \geq 1$ then the furthest point is just $\mathbf{d}$, and $r(\mathbf{d})$ is just $d_1 + \cdots + d_n - n + 1$, and otherwise we need minor corrections for those $d_i \leq 0$. Now we claim that $r$ is a cubism.

To show that $r$ attains only positive integer values, we can write $r$ as

$$r(\mathbf{d}) = 1 + \sum_{i=1}^{n} \max(d_i - 1, -d_i);$$

since $\max(d_i - 1, -d_i)$ is non-negative for any $d_i \in \mathbb{Z}$, $r$ attains only positive values. We leave the verification of (1) and (2) in the definition of a cubism to the reader. □

We also remark that—unlike the above example—there is no need for $r^{-1}(\{m\})$ to be finite; in fact, the next example shows that it can be convenient for $r^{-1}\{m\}$ to be infinite.

### 7.5. Other Examples of Cubisms and the Proof of Theorem 3.4.

*Proof of Theorem 3.4.* Let

$$\mathcal{D} = \{\mathbf{d} \mid a \leq \deg(\mathbf{d}) \leq a + n - 1\}.$$

Define $r \colon \mathbb{Z}^n \to \mathbb{N}$ as

$$r(\mathbf{d}) = \begin{cases} \deg(\mathbf{d}) - a + n + 1 & \text{if } \deg(\mathbf{d}) \geq a + n, \text{ and} \\ a + n - \deg(\mathbf{d}) & \text{if } \deg(\mathbf{d}) < a + n. \end{cases}$$

Setting $\mathcal{D}_0 = \mathcal{D}$ and, for $m \in \mathbb{N}$, $\mathcal{D}_m$ as in (63), we easily see that that if $r(\mathbf{d}) = m$ then $\mathrm{Cube}(\mathbf{d}) \setminus \mathcal{D}_{m-1}$ consists of a single point, namely $\mathbf{d}$ if $\deg(\mathbf{d}) \geq a + n$, and otherwise the single point $\mathbf{d} - \mathbf{1}$. We easily see that these single points are distinct as $\mathbf{d}$ varies over all $\mathbf{d} \notin \mathcal{D}$, and it follows that $r$ is a cubism of $\mathcal{D}$. □

**Example 7.7.** One can show by a cubism argument that the set $\mathcal{D} \subset \mathbb{Z}^2$ given by

$$\{(0,0)\} \cup \{\mathbf{d} \in \mathbb{Z}^2 \mid \deg(\mathbf{d}) = \pm 1\}$$

is fundamental, by defining $r(\mathbf{d})$ to be $|d_1|$ if $\deg(\mathbf{d}) = 1$ and otherwise $||\deg(\mathbf{d})|-1|$; we depict this cubism in Figure 3. It follows that any subset of $\mathcal{D}$ is subfundamental (e.g., removing $(0,0)$), and any superset of $\mathcal{D}$ is superfundamental.

It is intriguing—but not relevant to this article—to consider the various other fundamental modular domains of $\mathbb{Z}^n$.

We also note that in Example 7.7, it may be simpler to first extend a function $\mathcal{D} \to \mathbb{Z}$ along all points of degree 0, whereupon the extension is defined on all points of degree between $-1$ and $1$, and then further extend the function to all of $\mathbb{Z}^n$. In this case one can view the set of 2-cubes as a well-ordered set, where all points of degree 0 are ordered before all points of degrees not between $-1$ and $1$. One can therefore define a more general cubism as any well-ordering of the $n$-cubes of $\mathbb{Z}^n$, or, more generally, any partial ordering such that each subset of $n$-cubes has a minimal element. The proofs of all theorems easily generalize to these more general notions of a cubism.

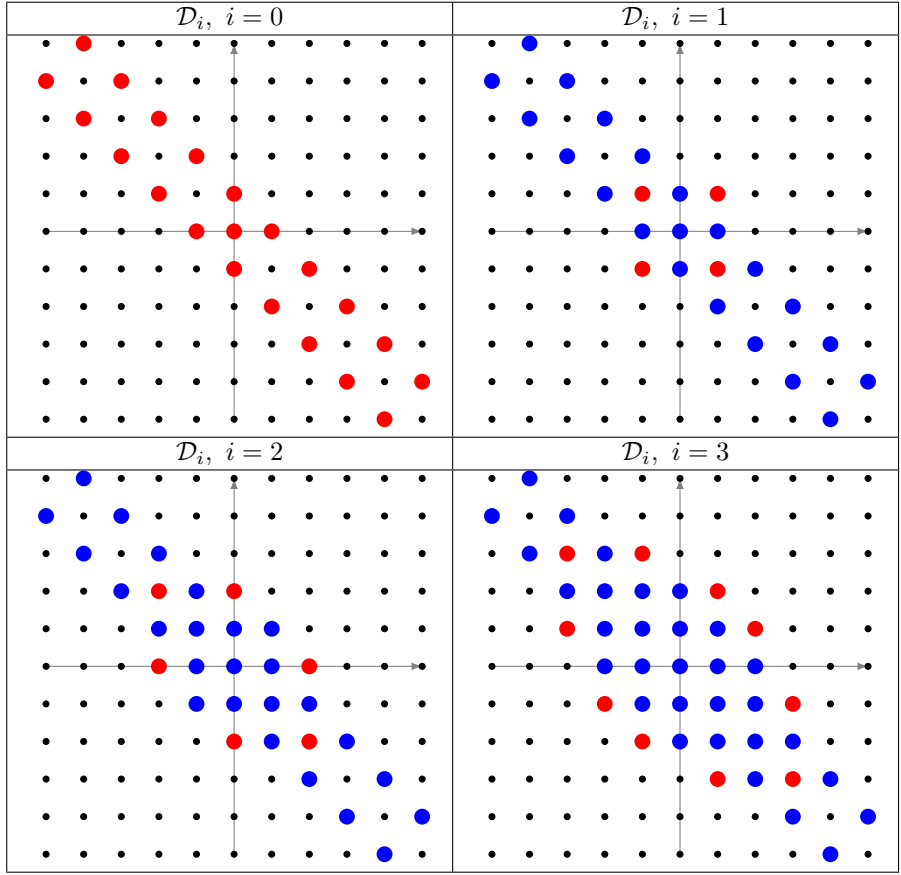

(A) New points in red, old points in blue

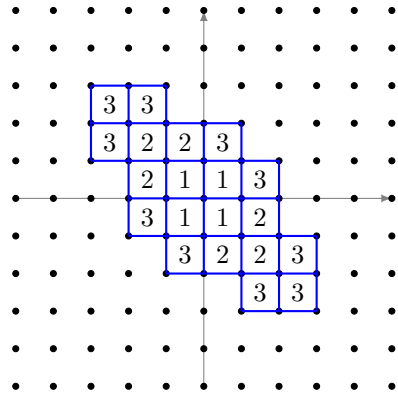

(B) The cubism after 4 steps.

FIGURE 3. A Cubism for Example 7.7.

Department of Mathematics, University of British Columbia, Vancouver, BC V6T 1Z2, CANADA.

*Email address*: nicolasfolinsbee@gmail.com

Department of Computer Science, University of British Columbia, Vancouver, BC V6T 1Z4, CANADA.

*Email address*: jf@cs.ubc.ca